\documentclass{article}
\usepackage{amsmath}
\usepackage{amsthm}
\usepackage{amsfonts}
\usepackage{amssymb}
\usepackage{latexsym}
\usepackage[dvips]{graphicx}
\usepackage{pstricks}
\usepackage{epsfig}
\usepackage{a4}
\usepackage{setspace}
\usepackage{color}
\usepackage{amsfonts,amssymb,graphics,epsfig,verbatim,bm,latexsym,amsmath,url,amsbsy}
\usepackage{amsfonts,amssymb,graphics,epsfig,verbatim,bm,latexsym,amsmath,url,amsbsy}

 \oddsidemargin=\evensidemargin
\newtheorem{theorem}{Theorem}

\newtheorem{proposition}{Proposition}

\newtheorem{remark}{Remark}

\begin{document}

\title{State-Dependent Autoregressive Models: Properties,
Estimation and Forecasting}
\author{Fabio Gobbi\footnote{Department of Statistics, University of Bologna, Italy} \quad
Sabrina Mulinacci\footnote{Department of Statistics, University of Bologna, Italy}}
\maketitle

\abstract{This paper studies some temporal dependence properties and addresses the issue of parametric estimation for a class of state-dependent autoregressive
models for nonlinear time series in which we assume a stochastic autoregressive coefficient depending on the first lagged value of the process itself. We call such a model
\emph{state-dependent first-order autoregressive process, (SDAR)}. We introduce some assumptions under which this class of models is strictly stationary and uniformly ergodic and we establish consistency and asymptotic normality of the quasi-maximum likelihood estimator of the parameters. In order to capture the potentiality of the model, we present an empirical application to nonlinear time series provided by the weekly realized volatility extracted from returns of some European financial indices. The comparison of forecasting accuracy is made considering an alternative approach provided by a two-regime SETAR model.}

Mathematics Subject Classification (2010): 60G10, 62M10, 91B84\\

JEL classification: C5, C01, C22, C58\\

{\bf Keywords}: Nonlinear time series, state-dependent autoregressive models, stationarity, ergodicity, quasi-maximum likelihood, forecasting accuracy.

\medskip

\section{Introduction}
In this paper we propose a generalized version of a first-order autoregressive process to model nonlinear time series where the autoregressive coefficient depends on the first
lagged state variable \begin{equation}\label{modelIntro}
Y_t=\alpha+\psi(Y_{t-1};\boldsymbol{\gamma}) Y_{t-1}+\xi_t,
\end{equation}
where $\psi$ is a specified function satisfying some assumptions and depending on a set of parameters $\boldsymbol{\gamma}$ and the error term $\xi_t$ is independent of $Y_{t-1}$ with zero mean and volatility $\sigma$.
The model is related to a much wider class of models with stochastic coefficients given by the family of \emph{functional autoregressive processes} of the form
$Y_t=f(Y_{t-1})+\xi_t$ in which the function $f(\cdot)$ can satisfy the most different properties. In this research line we have a number of contributions starting
from Hastie and Tibshirani (1990). Hardle et al. (1997) reviews alternative nonparametric estimation techniques and Diaconis and Feedman (1999) offer a method for
studying the steady state distribution of a Markov chain using iterated random functions and give useful bounds on rates of convergence in a variety of examples.
Finally, it worths to mention the two very exhaustive books on nonlinear time series of Fan and Yao (2003) and Douc et al. (2014).

In line with the functional autoregressive approach, Chen and Tsay (1993) explored the particular case of the "functional coefficient AR" model
$$Y_t=f_1(Y_{t-d})Y_{t-1}+\cdots +f_p(Y_{t-d})Y_{t-p}+\xi_t$$
with $d>0$ is some specified delay, which generlizes the well known ``exponential autoregressive (EXPAR)'' model introduced in Haggan and Ozaki (1981).
Cai et al (2000) adopt local linear regression techniques to estimate functional coefficient regression models for times series data.
Chen and Liu (2001) study nonparametric estimation and hypothesis testing procedures for the same model.

In this work we study the particular case, defined in (1.1), of the class of functional-coefficient AR models where $d=p=1$.
We call our class of models \emph{state-dependent} first-order autoregressive process, SDAR.
A preliminary version of this type of models can be found in Cherubini and Gobbi (2013) and in Cherubini et al. (2016) where the model is applied to time series of interest rates.
Our aim is to study the \emph{persistence function} $\psi(\cdot)$ in order to derive the necessary assumptions for the SDAR model to generate nonlinear time series that satisfy some serial dependence properties, such as  strict stationarity and uniform ergodicity, and such that the model parameters can be estimated using the quasi-maximum likelihood (QML) estimation technique. More precisely, we introduce some assumptions on $\psi$ in order to ensure that the QML estimator is consistent and asymptotically normal. From this point of view, the biggest problem is to assess that some functionals
involving the data process $Y_t$ and the vector of parameters $\boldsymbol{\theta}=(\alpha, \boldsymbol{\gamma}, \sigma)$ satisfy a strong unform law of large numbers.
For this aim, we use a result of Potscher and Prucha (1989) who proved a generic strong uniform law of large numbers for stochastic processes under general properties of serial dependence and heterogeneity.
Given the structure of the SDAR model, we apply the conditions needed in the theorem of Potscher and Prucha (other possible approaches can be found in
 Bierens, 1981 and 1984, and in Andrews, 1987).

Notice that the potential of the SDAR model is to be found above all in the fact that, having a stochastic autoregressive coefficient,
it generates nonlinear time series which are adequate to explain many characteristics observed in contemporary time series, such as
volatility clustering and extreme value dependence and moreover, it can be used to model temporal dependencies ranging from short to long-memory.

Moreover, as we will show, model (1.1) is a particular specification of the wider class of convolution-based autoregressive processes introduced in Cherubini et al (2016) and studied also in Gobbi and Mulinacci (2019).
Such models, considered in full generality, are untractable: aim of this paper is to consider a specific case in which temporal dependence properties and
estimation issues can be established.

In order to appreciate the potentiality of proposed SDAR models we present an empirical application to real economic nonlinear time series. In particular, we consider the weekly realized volatility extracted from returns of three European financial indices, CAC40 (France), DAX30 (Germany) and FTSE100 (UK). We formulate and estimate a SDAR model on a sub-sample of the historical data and the remaining out-of sample data are considered for a forecasting perspective. In particular, we compare forecast accuracy of the selected SDAR models with an alternative approach for modelling nonlinear time series given by the self-exciting threshold autoregressive models (SETAR), which were first proposed and studied by Tong (1978, 1986 and 1995) and Tong and Lim (1980). In SETAR models the variable $y_t$ is a linear autoregression within a regime but may move among regimes depending on the value taken by a lag of $y_t$ itself. SETAR models have been applied to a number economic and financial variables. For example, and among others, Krager and Kugler (1993), Peel and Speight (1994) and Chappell et al. (1996) apply such models to exchange rates and Tiao and Tsay (1994), Potter (1995) and Clements and Krozling (1998) to U.S. GNP.

The plan of the paper is the following. Section 2 introduces SDAR models and discusses their theoretical properties such as stationarity and ergodicity and the method of estimation. In section 3 we present an empirical application to realized volatility and we compare forecasting accuracy of SDAR models with leading nonlinear time series models such as the SETAR model. Section 4 concludes.

\section{SDAR models: definitions and theoretical results}
In this section we introduce the SDAR model and we present its most significant theoretical results such as strictly stationarity, uniform ergodicity and the consistency and the asymptotic normality of the QML estimator of the vector of parameters.\\
Let $(Y_t)_{t \in \mathbb{N}}$ be a stochastic process defined on a complete probability space
$(\mathbb{R}^{\infty}, \mathfrak{B}(\mathbb{R}^{\infty}), \mathbb{P}^0)$. We say that $(Y_t)_t$ is a SDAR process if it satisfies the following
specification

\begin{equation}\label{model}
\left\{\begin{array}{cc}
                                       Y_{t}=\alpha+\psi(Y_{t-1};\boldsymbol{\gamma})Y_{t-1}+\xi_t,& t \geq 1,\\
                                      \xi_t \sim i.i.d. \ N(0,\sigma),\\
                                      \end{array}\right ..
\end{equation}

where $Y_0 = 0$ and $\psi(Y_{t-1};\boldsymbol{\gamma})$
is a measurable function of the lagged variable $Y_{t-1}$ that depends on a $p$-dimensional vector of
parameters $\boldsymbol{\gamma}$. $\psi$ specifies a dynamics for the autoregressive coefficient which is not longer a
constant as in the standard AR(1) case: we call this function \emph{persistence function}. The sequence
of error terms, $(\xi_t)_t$, is i.i.d. and normally distributed with zero mean and volatility $\sigma$.\\ The model is a specific case of the class of convolution-based autoregressive processes (see Cherubini et al., 2016) defined as
$$Y_t=\alpha+\phi Y_{t-1}+\eta_t$$
where $\eta_t$ \emph{depends} on $Y_{t-1}$. More precisely, it corresponds to the case in which
$$\eta_t=h(Y_{t-1})+\xi_t$$
with $h=\psi(y)y-\phi y$ and $(\xi_t)_t$ i.i.d. and normally distributed $N(0,\sigma)$.\bigskip

We denote by $\boldsymbol{\theta}=(\alpha,\boldsymbol{\gamma},\sigma)$ the vector of model parameters belonging to
$\Theta=\mathcal{A} \otimes \Gamma \otimes \mathcal{S} \subset \mathbb{R}^{p+2}$ and we set $Y^t=Y_t|Y_{t-1}$.
Moreover, let $p_t^{\boldsymbol{\theta}}(Y^t)$ be the conditional density of $Y_t$ given $Y_{t-1}$ and let $\boldsymbol{\theta}^0$ be the true
value of the parameter. The assumptions required case by case by the main results, are the following:
\begin{itemize}
\item \textbf{a1}. $\psi$ is differentiable with respect to $y$ and
$|\psi(y;\boldsymbol{\gamma})|+\left|y\frac{d}{dy}\psi(y;\boldsymbol{\gamma})\right|\leq K<1$ $\forall y \in \mathbb{R}$
  \item \textbf{a2}. $\psi(y;\boldsymbol{\gamma})y$ is uniformly bounded in $y$.

\end{itemize}
\begin{remark} Notice that assumption \textbf{a1} implies that the persistence function has to assume values far away from 1.\end{remark}

\begin{remark}\label{moments}
Under assumption {\bf a2}, all absolute moments of $Y_t$ of any order are finite.
\end{remark}
\vspace{0.2cm}
We shall now introduce two propositions which establish that the SDAR model is strictly stationary and uniformly ergodic for any ${\boldsymbol\theta}\in \Theta$. Such properties are desirable in light
of the empirical potential of this family of models. Proofs are reported in the appendix 1.

\begin{proposition}\label{stationarity}
Under assumption  {\bf a1}, the SDAR model in (2.2) is strictly stationary.
\end{proposition}
\begin{proof} See Appendix 1.
\end{proof}

\begin{proposition}\label{ergodicity}
Under assumption {\bf a2}, the SDAR model in (2.2) is uniformly ergodic.
\end{proposition}
\begin{proof} See Appendix 1.
\end{proof}

As regards the estimation methodology of the vector of parameters $\boldsymbol{\theta}=(\alpha, \boldsymbol{\gamma}, \sigma)$ which characterize SDAR models we are going to establish the consistency and the asymptotic normality of the QML estimator.
In order to make more easier the reading the statement of the main theorem and the proofs in the appendix 1,
we state the notation.\\
From an estimation point of view, we are interested in partial derivatives with respect to the parameters. The gradient of the function $\psi$ is denoted by
$\nabla_{\boldsymbol{\gamma}} \psi(Y_{t-1};\boldsymbol{\gamma})=\left[ \frac{\partial}{\partial \gamma_k}  \psi(Y_{t-1};\boldsymbol{\gamma}) \right]_{k=1,...,p}=\left[ \psi^{\gamma_k}(Y_{t-1};\boldsymbol{\gamma}) \right]_{k=1,...,p}.$
The hessian matrix of $\psi$ is
$$\nabla_{\boldsymbol{\gamma}}^2 \psi(Y_{t-1};\boldsymbol{\gamma})=\left[ \frac{\partial^2}{\partial \gamma_k \partial \gamma_j}  \psi(Y_{t-1};\boldsymbol{\gamma}) \right]_{k,j=1,...,p}=\left[ \psi^{\gamma_k,\gamma_j}(Y_{t-1};\boldsymbol{\gamma}) \right]_{k,j=1,...,p}.$$
Furthermore, given a time series $\boldsymbol{y}^n=(y_1,...,y_n)$ generated by the SDAR model in (2.2) the quasi log-likelihood associated to $\boldsymbol{y}^n$ is
$$L_n(\boldsymbol{y}^n;\boldsymbol{\theta})=\frac{1}{n}\sum_{t=1}^n \ell_t(Y^t;\boldsymbol{\theta}),$$ where $\ell_t(Y^t;\boldsymbol{\theta})=\ln p_t^{\boldsymbol{\theta}}(Y^t)$. More explicitly,
since the conditional distribution of $Y_t$ given $Y_{t-1}$ is $N\left(\alpha+\psi(Y_{t-1};\boldsymbol{\gamma})Y_{t-1},\sigma\right)$ the conditional density $p_t^{\boldsymbol{\theta}}(Y^t)=\frac{1}{\sigma \sqrt{2 \pi}}e^{-\frac{(Y_t-\alpha-\psi(Y_{t-1};\boldsymbol{\gamma})Y_{t-1})^2}{2\sigma^2}}$ and hence
$$\ell_t(Y^t;\boldsymbol{\theta})= \ln p_t^{\boldsymbol{\theta}}(Y^t)
\propto -log(\sigma)-\frac{(Y_t-\alpha-\psi(Y_{t-1};\boldsymbol{\gamma})Y_{t-1})^2}{2 \sigma^2}.$$
Hence, the gradient and the hessian matrix of $\ell_t(Y^t;\boldsymbol{\theta})$ are respectively
$$\underset{(p+2) \times 1}{\nabla_{\boldsymbol{\theta}}\ell_t(Y^t;\boldsymbol{\theta})}
=\left(
 \begin{array}{c}
   \ell_t^{\alpha}(Y^t;\boldsymbol{\theta})\\
  \underset{(p \times 1)}{\nabla_{\boldsymbol{\gamma}}\ell_t(Y^t;\boldsymbol{\theta})} \\
  \ell_t^{\sigma}(Y^t;\boldsymbol{\theta}) \\                                                                              \end{array}
     \right)
$$
where $\ell_t^{\alpha}(Y^t;\boldsymbol{\theta})=\frac{\partial}{\partial \alpha}  \ell_t(Y^t;\boldsymbol{\theta})$,
$\ell_t^{\sigma}(Y^t;\boldsymbol{\theta})=\frac{\partial}{\partial \sigma}  \ell_t(Y^t;\boldsymbol{\theta})$ and
$\nabla_{\boldsymbol{\gamma}}$ denotes the gradient with respect to $\boldsymbol{\gamma}$, and

$$\underset{(p+2) \times (p+2)}{\nabla^{2}_{\boldsymbol{\theta}}\ell_t(Y^t;\boldsymbol{\theta})}=
$$$$=\left(
   \begin{array}{ccccc}
\ell_t^{\alpha \alpha}(Y^t;\boldsymbol{\theta}) & \ell_t^{\alpha \gamma_1}(Y^t;\boldsymbol{\theta})
& \ldots & \ell_t^{\alpha \gamma_p}(Y^t;\boldsymbol{\theta}) &  \ell_t^{\alpha \sigma}(Y^t;\boldsymbol{\theta})\\
\ell_t^{\gamma_1 \alpha}(Y^t;\boldsymbol{\theta}) & \ldots
& \underset{p \times p}{\nabla_{\boldsymbol{\gamma}}^2 \ell_t(Y^t;\boldsymbol{\theta})}& \ldots &  \ell_t^{\gamma_1 \sigma}(Y^t;\boldsymbol{\theta})\\
\vdots & \ddots & \vdots & \ddots & \vdots\\
\ell_t^{\sigma \alpha}(Y^t;\boldsymbol{\theta}) & \ell_t^{\sigma \gamma_1}(Y^t;\boldsymbol{\theta})
& \ldots & \ell_t^{\sigma \gamma_p}(Y^t;\boldsymbol{\theta}) &  \ell_t^{\sigma \sigma}(Y^t;\boldsymbol{\theta})\\
   \end{array}
 \right).
$$
The QML estimator
is the solution of the maximization problem
$$\hat{\boldsymbol{\theta}}_n=\arg\max_{\boldsymbol{\theta} \in \Theta} \frac{1}{n}\sum_{t=1}^n \ell_t(Y^t;\boldsymbol{\theta}).$$
Consistency and asymptotic normality of this estimator require three additional assumptions.
\begin{itemize}
  \item \textbf{a3}. The parameter space
$\Theta$ is a compact subset of $\mathbb{R}^{p+2}$.
  \item \textbf{a4}.
First and second-order partial derivatives with respect to the parameters of the persistence function $\psi$ are continuous and
uniformly bounded in $y$.
  \item \textbf{a5}.
$\vert \psi^{\gamma _k}(y;\boldsymbol{\gamma})y\vert\leq C$ uniformly on $\mathbb R\times\Theta$, for all $k$ and
$\vert \psi^{\gamma _k\gamma_j}(y;\boldsymbol{\gamma})y\vert\leq D$ uniformly on $\mathbb R\times\Theta$, for all $k,j$.
\end{itemize}
Let us introduce the main theorem.
\begin{theorem}\label{QMLE}
Under assumptions \textbf{a1}-\textbf{a5} the estimator
$\hat{\boldsymbol{\theta}}_n$ is strongly consistent for $\boldsymbol{\theta}^0$ and moreover it satisfies
$$\bar{H}^0_{n} \left(G^0_{n}\right)^{-1/2}\sqrt{n}(\hat{\boldsymbol{\theta}}_n-\boldsymbol{\theta}^0)
\stackrel{d}{\longrightarrow} N(0,I),$$
where $$\bar{H}^0_{n}=\mathbb{E}\left[\frac{1}{n}\sum_{t}\nabla^{2}_{\boldsymbol{\theta}}\ell_t(Y^t;\boldsymbol{\theta}^0)\right],$$ $$G^0_{n}=\mathbb{E}\left[\frac{1}{n}\sum_{t}\left(\nabla_{\boldsymbol{\theta}}\ell_t(Y^t;\boldsymbol{\theta}^0)\right)
\left(\nabla_{\boldsymbol{\theta}}\ell_t(Y^t;\boldsymbol{\theta}^0)\right)^T
\right]$$ and $I$ is the identity matrix of order $(p+2)$.
\end{theorem}
\begin{proof} See Appendix 1.
\end{proof}

\section{Empirical application}
In this section we estimate the SDAR model using empirical time series by selecting two different functional form of the persistence function $\psi$. The forecasting performance is evaluated using SETAR models as the benchmark.

\subsection{The data set}
Our application considers the weekly realized volatility extracted from daily returns of three different European financial market indices: CAC40 (France), DAX30 (Germany) and FTSE100 (U.K.). The sample period goes from January 2004 until April 2019. We consider a time series $\boldsymbol{r}=(r_1,...,r_N)$ of $N=3890$ daily returns from 2004.1 to 2018.48 for estimation purposes whereas the remaining observations from 2018.49 to 2019.16 will be used for forecasting perspective. We recover the corresponding in-sample time series of weekly realized volatility considering the square root of the sum of squared returns within each week, $v_t=\sqrt{\sum_{s=5(t-1)+1}^{5t} r^2_{s}}$, $t=1,...n$ with $n=778$. We apply the SDAR model to the logarithm of the volatility, $y_t=\ln(v_t)$.

We make use of various tests to asses if time series of log weekly realized volatility are nonlinear. In literature there are a number of nonlinearity tests which can be applied, but we concentrate on four of them. The first two are used for testing the neglected nonlinearity in the case where the null is the hypotheses of linearity in mean: the Teraesvirta Neural Network test (tnn-test) introduced in Teraesvirta et al. (1993) and the White Neural Network test (wnn-test) discussed in Lee et al. (1993). The third one is the likelihood ratio test for threshold nonlinearity (tlrt-test) implemented by Chan (1990). The null hypothesis is that the fitted model to the time series is an AR model with a specified lag structure and the alternative is that the fitted model is a threshold autoregressive model with the same lag structure for each regime.

 The last test we implement is a test for quadratic nonlinearity in a time series in which the null hypothesis is a normal AR process. The test (tsay-test) was introduced and implemented in Tsay (1986). We use R packages "fNonlinear" and "TSA" to perform the four tests. The results are presented in table 1.
 They are interesting and encouraging at the same time in the sense that we may reject the null in a number of cases
 but we observe strong evidence of nonlinear components only when we consider a lag equal to 1. For lag structures of higher order only the realized volatility extracted from the FTSE100 highlights strong evidence of nonlinearity whereas in the case of CAC40 and DAX30 such evidence is not found in particular considering tnn-test and wnn-test.
 This may have consequences on the predictive ability of the models used since the degree of nonlinearity present in the time series affects, how easy to infer, the performance of nonlinear models.

\begin{table}
\begin{center}
\begin{tabular}{cccc}
\hline
&   CAC40    &    DAX30    &   FTSE100  \\
\hline
&  &  lag=1  &  \\
\hline
tnn-test   & 0.0004    & 0.0019  &   $<$10$^{-5}$ \\
wnn-test   &  0.0004   & 0.0006  &$<$10$^{-5}$\\
tlrt-test  &  0.0009   &  0.0036   &   $<$10$^{-5}$\\
tsay-test  &   0.0001   &  0.0004   &  $<$10$^{-5}$  \\
\hline
&  &  lag=2  &  \\
\hline
tnn-test   & 0.1486    & 0.0172  &   0.0016 \\
wnn-test   &  0.0845   & 0.2886  &   0.0073\\
tlrt-test  &  0.0257   &  0.0028   &  0.0016\\
tsay-test  &  0.0288    &  0.0023   &  0.0002      \\
\hline
&  &  lag=3 &  \\
\hline
tnn-test   & 0.4528    & 0.2408  &    0.0563\\
wnn-test   &  0.1228   & 0.9287  &   0.0105\\
tlrt-test  &  0.0337   &  0.1166   &   0.0083\\
tsay-test  &  0.1055    & 0.0311 &     0.0049          \\
\hline
\end{tabular}
\caption{Nonlinearity tests: $p$-values for different lag structures.}\label{table1}
\end{center}
\end{table}
\bigskip

\subsection{Two different choices for the persistence function $\psi$}
We consider two possible specifications of the persistence function, say $\psi_1$ and $\psi_2$,
satisfying assumptions \textbf{a1}-\textbf{a5}. We refer to the SDAR model characterized by $\psi_i$ as SDAR(Mi) with $i=1,2$. Below we introduce the functional forms of $\psi_1$ and $\psi_2$ and we discuss conditions under which assumptions \textbf{a1}-\textbf{a5} are satisfied.

\begin{itemize} \item \textbf{SDAR(M1)}. Consider the persistence function
$$\psi_1(y;\boldsymbol{\gamma})=e^{-(\gamma_0+\gamma_1 y^{2r})}, \quad  \gamma_0\in\mathbb R,\gamma_1, r >0,$$
characterized by the vector of parameters $\boldsymbol{\gamma}=(\gamma_0, \gamma_1, r)$.
Under this specification, the SDAR model is a special case of the EXPAR model in Chen and Tsay (1993). Notice that if $\gamma_1=0$ we recover an AR(1) model.
\\
It can be easily verified that requirements {\bf a1} and {\bf a2} are satisfied since
$$\underset{y\in\mathbb R}\sup\left (|\psi_1(y;\boldsymbol{\gamma})|+\left|y\frac{d}{dy}\psi_1(y;\boldsymbol{\gamma})\right|\right )=
2re^{-\frac{2r\gamma_0+2r-1}{2r}}.$$
{\bf a4} and {\bf a5} are satisfied on any compact set
$$\Theta\subset\left \{(\gamma_0,\gamma_1)\in (\ln(2r)-1+\frac 1{2r},+\infty)\times(0,+\infty)\right\}$$
since the partial derivatives with respect to the parameters are the following: $\psi^{\gamma_0}_1(y;\boldsymbol{\gamma})=-\psi_1(y;\boldsymbol{\gamma})$,
$\psi^{\gamma_1}_1(y;\boldsymbol{\gamma})=-y^{2r}\psi_1(y;\boldsymbol{\gamma})$,
$\psi^{r}_1(y;\boldsymbol{\gamma})=-\gamma_1 y^{2r} \ln(y^2)\psi_1(y;\boldsymbol{\gamma})$
$\psi^{\gamma_0 \gamma_0}_1(y;\boldsymbol{\gamma})=\psi_1(y;\boldsymbol{\gamma})$,
$\psi^{\gamma_0 \gamma_1}_1(y;\boldsymbol{\gamma})=y^{2r} \psi_1(y;\boldsymbol{\gamma})$,
$\psi^{\gamma_0 r}_1(y;\boldsymbol{\gamma})=\gamma_1 y^{2r} \ln(y^2) \psi_1(y;\boldsymbol{\gamma})$,
$\psi^{\gamma_1 \gamma_1}_1(y;\boldsymbol{\gamma})=y^{4r}\psi_1(y;\boldsymbol{\gamma})$,
$\psi^{\gamma_1 r}_1(y;\boldsymbol{\gamma})=y^{2r} \ln(y^2) \psi_1(y;\boldsymbol{\gamma})\left[\gamma_1 y^{2r}-1\right]$,
$\psi^{r r}_1(y;\boldsymbol{\gamma})=\gamma_1 y^{2r} \ln^2(y^2) \psi_1(y;\boldsymbol{\gamma})\left[\gamma_1 y^{2r}-1\right]$.
%Since the persistence function $\psi_1$ depends on $y$ the interpretation of the coefficients $\gamma_0$ and $\gamma_1$ can be rather hard in order to establish their role in the construction of the process and therefore in the temporal dependence that characterize the time series. Nevertheless, we can observe that $\gamma_0$ is a parameter that influence the absolute value of the persistence function in the sense that if $\gamma_0$ is close to zero the value of the persistence function reaches one if $y$ is sufficiently close to zero. The speed with which this happens depends on the value of $\gamma_1$. Clearly, this affects short and long-memory and mean-reverting properties of the time series. On the other hand, $\gamma_1$ is a shape parameter, i.e., $\psi_1$ decreases faster when $y$ moves away from zero as $\gamma_1$ grows and this mainly affects the memory of the process.

%We try to describe the dynamics of the persistence functions $\psi_1(y;\boldsymbol{\gamma})$ as a function of $y$ setting the shape parameters very close to those estimated in the following sections, namely, $\gamma_0=0.4$, $\gamma_1=0.07$ and $r=0.3$, so as to have a clear perception of how the time-dependent autoregressive coefficient follows the volatility trend. Despite the range of possible values is rather narrow (approximately 0.60-0.70 ), the dynamics is monotone increasing with respect to the volatility. In other words, periods (weeks in our case) of high volatility will be characterized by increasing values of $\psi_1(y;\boldsymbol{\gamma})$ indicating greater persistence.

\item  \textbf{SDAR(M2)}. Consider the persistence function
  $$\psi_2(y;\boldsymbol{\gamma})=\frac{1}{\gamma_0+\gamma_1 y^{2r}}, \quad  \gamma_0>1, \ \gamma_1 >0.$$
Requirements {\bf a1} and {\bf a2} are satisfied since
$$\underset{y\in\mathbb R}\sup\left (|\psi_1(y;\boldsymbol{\gamma})|+\left|y\frac{d}{dy}\psi_1(y;\boldsymbol{\gamma})\right|\right )=
\frac{(1+2r)^2}{8r\gamma_0}.$$
{\bf a4} and {\bf a5} are satisfied on any compact set
$$\Theta\subset\left \{(\gamma_0,\gamma_1)\in \left(\frac{(1+2r)^2}{8r},+\infty\right)\times\left(0,+\infty\right)\right\}$$ since the partial derivatives with respect to the parameters are the following:
 $\psi^{\gamma_0}_2(y;\boldsymbol{\gamma})=-\psi^{2}_2(y;\boldsymbol{\gamma})$,
$\psi^{\gamma_1}_1(y;\boldsymbol{\gamma})=-y^{2r} -\psi^{2}_2(y;\boldsymbol{\gamma})$,
$\psi^{r}_1(y;\boldsymbol{\gamma})=-\gamma_1 y^{2r} \ln(y^2) \psi^{2}_2(y;\boldsymbol{\gamma})$,
$\psi^{\gamma_0 \gamma_0}_2(y;\boldsymbol{\gamma})=2 \psi^{3}_2(y;\boldsymbol{\gamma})$,
$\psi^{\gamma_0 \gamma_1}_2(y;\boldsymbol{\gamma})=2 y^{2r} \psi^{3}_2(y;\boldsymbol{\gamma})$,
$\psi^{\gamma_0 r}_2(y;\boldsymbol{\gamma})=2 \gamma_1 y^{2r} \ln(y^2)\psi^{3}_2(y;\boldsymbol{\gamma})$,
$\psi^{\gamma_1,\gamma_1}_2(y;\boldsymbol{\gamma})=y^{4r} \psi^{3}_2(y;\boldsymbol{\gamma})$,
$\psi^{\gamma_1 r}_2(y;\boldsymbol{\gamma})=y^{2r} \ln(y^2)\psi^{2}_2(y;\boldsymbol{\gamma})\left[2 \gamma_1\psi_2(y;\boldsymbol{\gamma})-1 \right]$,\\
$\psi^{r r}_2(y;\boldsymbol{\gamma})=\gamma_1 y^{2r} \ln^2(y^2)\psi^{2}_2(y;\boldsymbol{\gamma})\left[2 \gamma_1 y^{2r} \psi_2(y;\boldsymbol{\gamma})-1 \right]$.\\
 % The dynamics of $\psi_2(y;\gamma_0,\gamma_1,r)$ setting $\gamma_0=1.2$, $\gamma_1=0.08$ and $r=0.5$ is still monotone increasing but the range of possible values is wider than in the previous case showing a significantly higher persistence.
\end{itemize}

\subsection{Estimation results}
We present and discuss estimation results applying the QML method introduced in section 2. We estimate the parameters for both models SDAR(M1) and SDAR(M2) for the realized volatility extracted from CAC40, DAX30 and FTSE100. The parameters estimates are summarized in table 2. Some insights are possible observing these results. Firstly, and independently of the adopted model, all estimates are highly significant indicating that all parameters of the model are relevant. Moreover, as expected, the estimate of $r$ is systematic higher for the SDAR(M2) model with respect to the SDAR(M1) model. A possible explanation of this result depends on the fact that the specification characterized by $\psi_1$ is more sensitive to the variation of the parameter $r$. Furthermore, as it is easily to check, in order to reproduce the same value of the persistence (i.e., the same value of the coefficients $\psi_1$ and $\psi_2$) the parameter $r$ must be greater in the case M1 than in the case M2.

In order to identify the best model we use the AIC criterion, in the sense that for each financial index we compute the AIC associated to each model and we select the model corresponding to the minimum AIC. Table 3 summarize the results. With the AIC criterion we identify SDAR(M1) for CAC40 and DAX30 and SDAR(M2) for FTSE100. Furthermore, table 4 summarizes the evaluation of the goodness of the selected models through an independent test of the residuals. We use the BDS test developed by Brock et al. (1987) and later published by Brock et al. (1996). The test can be used as a portmanteau test of independent and identically distribution when applied to the residuals of a fitted model as in our case (see, among others, Krager and Kluger (1993)). In particular, if the null of i.i.d. residuals cannot be rejected, this indicates that the model provides an adequate fit to the original time series and successively removes the nonlinearity in the data. All $p$-values reported in table 4 allow to accept the null and therefore our selected SDAR models are adequate for all three indices.\\
It can be interesting to analyze the variation of the estimated time-varying persistence functions $\psi_1(y_{t-1}, \hat{\boldsymbol{\gamma}})$ and $\psi_2(y_{t-1}, \hat{\boldsymbol{\gamma}})$ along the sample period to investigate if they capture the volatility dynamics. Figures 1-3 confirm this expectation. In periods when weekly volatility is high, both functions increase indicating a greater persistence. This happens for all three indices. However for CAC40 the maximum value reached by the persistence function at the peak of volatility (which occurs in 2009, the year of the post-Lehman crisis for all three indices) is much lower than that of the other two indices DAX30 and FTSE100 indicating an excess persistence for the last two financial indices.

For forecasting purposes we estimate an alternative model widely used in nonlinear time series forecasting: the self-exciting autoregressive model (SETAR) introduced by Tong (1978) and intensively used for forecasting perspective in a number of articles as mentioned in the introduction of this paper. In appendix 2 we briefly describe the mathematical formulation of this model. The estimation procedure (Tong (1983, 1995), Potter (1995), among others) starts with the choice of the number of regimes $p$. We consider two cases $p=2$ or $p=3$. Once the number of regimes has been set, we choose a maximum value of the number of lags (for each regime) and proceed with the estimation of the models. The preferred model is that which minimizes AIC. We identify the following models:  SETAR(2,3,3) for the CAC40, SETAR(2,2,3) for the DAX30 and SETAR(2,4,4) for the FTSE100. Results concerning the
selected models are reported in the Appendix 2. As in the case of SDAR models we perform the BDS test for fitted residuals and the obtained $p$-values corroborate the goodness-of-fit of the estimated models.

\begin{table}
\begin{center}

\begin{tabular}{c|c|cc|cc}
\hline
       &    &   SDAR(M1) & & SDAR(M2) & \\
  \hline
 Index  & Parameters   &   Coeff. &  Std error  & Coeff. &  Std error    \\
  \hline

  CaC40  & $\alpha$   & -1.5856$^{***}$          &  0.0408          &  -1.3663$^{***}$        &  0.1218            \\
         & $\gamma_0$ & 0.3734$^{***}$          &   0.0157         &  1.1808$^{***}$        &  0.0189            \\
         & $\gamma_1$ &  0.0649$^{***}$         &   0.0049         &   0.0785$^{***}$       &  0.0067            \\
         & $\sigma$   &   0.5134$^{***}$        &   0.0192        &  0.5378$^{***}$        &   0.0112           \\
         & $r$        &   0.3198$^{***}$        &    0.0275        &  0.5596$^{***}$        &  0.0336            \\
  \hline
 DAX30  & $\alpha$   &  -1.8863$^{***}$       &   0.0512         &  -1.4375$^{***}$        & 0.0691  \\
         & $\gamma_0$ & 0.4453$^{***}$         &  0.0043          &  1.1346$^{***}$        &  0.0207   \\
         & $\gamma_1$ & 0.0736$^{***}$          &  0.0015          &  0.0973$^{***}$        & 0.0009    \\
         & $\sigma$   &  0.5134$^{***}$         &  0.0193          & 0.5524$^{***}$         &  0.0159      \\
         & $r$        &  0.4036$^{***}$         &   0.0155         &  0.5628$^{***}$        &  0.0525    \\
  \hline
FTSE100  & $\alpha$   & -1.8754$^{***}$          & 0.0574           &   -1.3753$^{***}$       & 0.0831    \\
         & $\gamma_0$ & 0.3701$^{***}$          &  0.0078          &   1.1705$^{***}$       &   0.0124    \\
         & $\gamma_1$ & 0.0945$^{***}$          &  0.0053          &   0.0884$^{***}$       & 0.0034      \\
         & $\sigma$   & 0.5041$^{***}$          &   0.0281         &   0.5017$^{***}$       &  0.0195     \\
         & $r$        &  0.3315$^{***}$         &   0.0284         &  0.4555$^{***}$        & 0.0498     \\
  \hline

\end{tabular}
\caption{Estimated parameters and relative standard errors of the adopted SDAR models for weekly realized volatilities obtained from the selected financial indices. Three asterisks denote that the parameter is significantly different from zero at the 1\% level.}\label{table2}
\end{center}
\end{table}
\bigskip

\begin{table}
\begin{center}

\begin{tabular}{c|c|c|c}
\hline
       &    &   SDAR(M1)  & SDAR(M2)  \\
  \hline
 Index    &   Fit  &   &        \\
  \hline

  CAC40     & AIC             &     1124.54        &   1134.30          \\

 \hline
    DAX30          & AIC              &     1148.82        &   1157.05         \\
 \hline

  FTSE100             & AIC              &     1151.36        &   1135.19          \\
 \hline

\end{tabular}
\caption{Measures of goodness of fit.}\label{table3}
\end{center}
\end{table}
\bigskip

\begin{table}
\begin{center}

\begin{tabular}{c|c|c|c}
  \hline
 BDS Test ($p$-values)   &  CAC40 &  DAX30  &  FTSE100   \\
  \hline

       &  eps[1] m=2: 0.7116   &  eps[1] m=2: 0.4390        &  eps[1] m=2: 0.1595    \\
            &  eps[1] m=3: 0.1581  &  eps[1] m=3: 0.2485     &  eps[1] m=3: 0.2838    \\
         &  eps[2] m=2: 0.4945  &  eps[2] m=2: 0.3058    &   eps[2] m=2: 0.1735      \\
          &  eps[2] m=3: 0.1735    &  eps[2] m=3: 0.1482   &   eps[2] m=3: 0.1957   \\
          &  eps[3] m=2: 0.4059    &  eps[3] m=2: 0.2029    &  eps[3] m=2: 0.1262    \\
            &  eps[3] m=3: 0.0855    &  eps[3] m=3: 0.0898   &  eps[3] m=3: 0.1385    \\
            &  eps[4] m=2: 0.3300      &  eps[4] m=2: 0.0895   &   eps[4] m=2: 0.1523      \\
            &  eps[4] m=3: 0.1150   &  eps[4] m=3: 0.0965    &  eps[4] m=3: 0.1670         \\
 \hline

\end{tabular}
\caption{BDS test for independence of residuals. This table reports the $p$-value output of the function bdsTest() of R. The bdsTest test examines the spatial dependence of the observed series. To do this, the series is embedded in $m$-space and the dependence is examined by counting near points. Points for which the distance is less than a constant EPS are called near. In this case:
embedding dimension $m=3$, EPS=(0.5*sd(residuals), sd(residuals), 1.5*sd(residuals), 2*sd(residuals)). Selected models: SDAR(M1) for CAC40 and DAX30, SDAR(M2) for FTSE100.}\label{table4}
\end{center}
\end{table}
\bigskip

\begin{figure}[htbp]
\centering
\includegraphics[width=10cm, height=6cm]{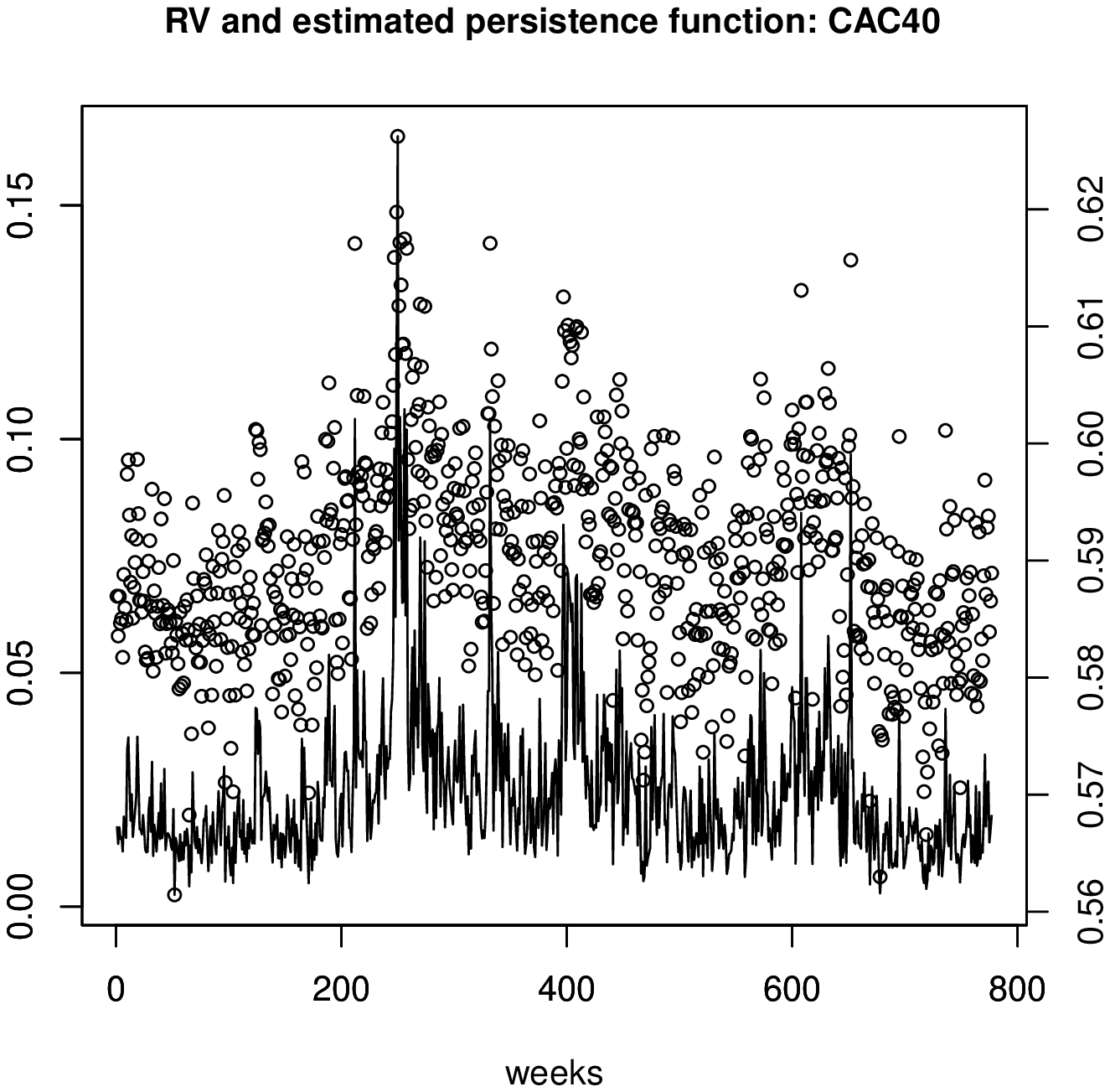}
\caption{Realized volatility (line, left vertical axis) and estimated persistence function (points, right vertical axis). SDAR(M1) for CAC40.}\label{fig1}
\end{figure}

\begin{figure}[htbp]
\centering
\includegraphics[width=10cm, height=6cm]{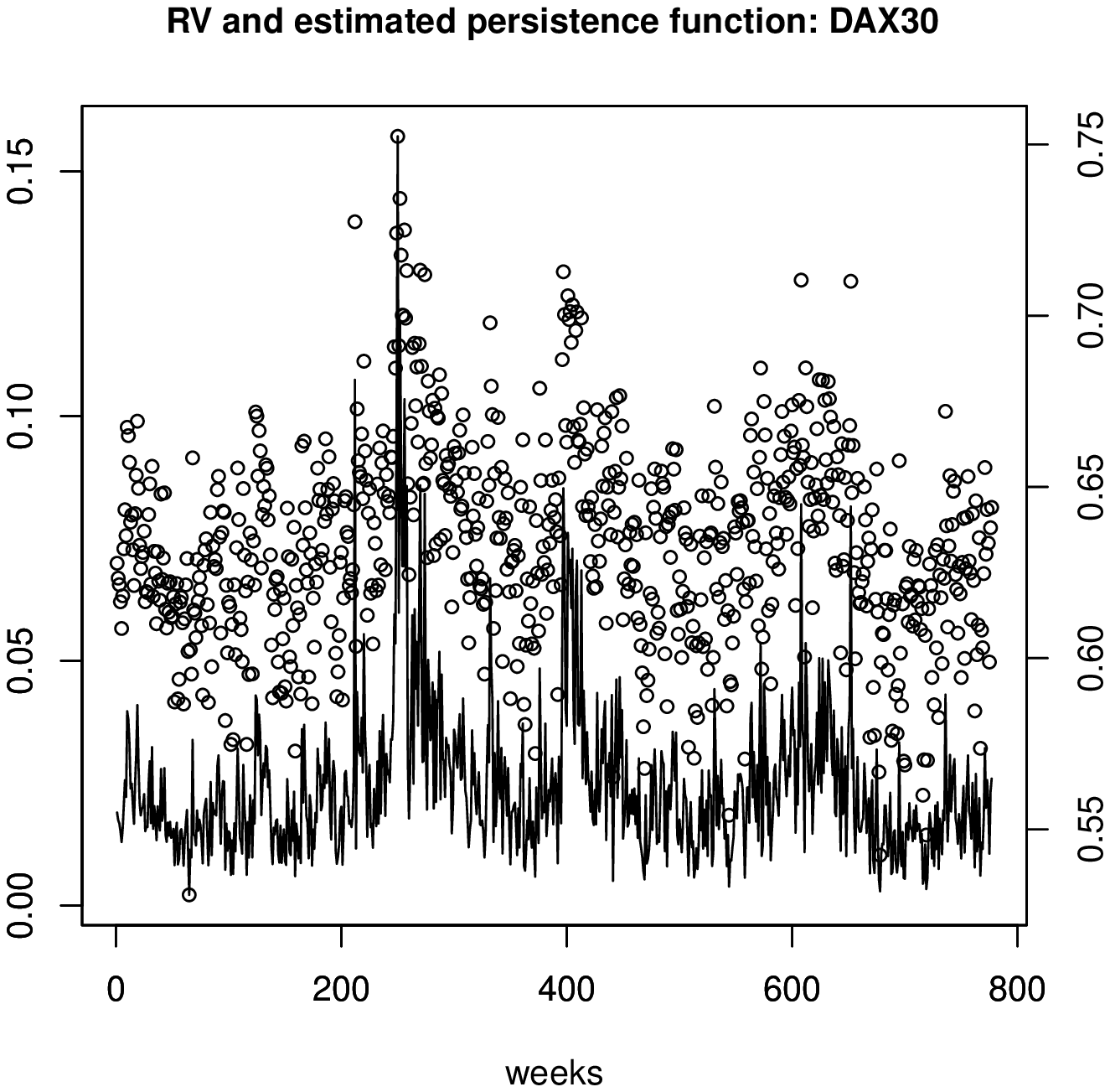}
\caption{Realized volatility (line, left vertical axis) and estimated persistence function (points, right vertical axis). SDAR(M1) for DAX30.}\label{fig2}
\end{figure}

\begin{figure}[htbp]
\centering
\includegraphics[width=10cm, height=6cm]{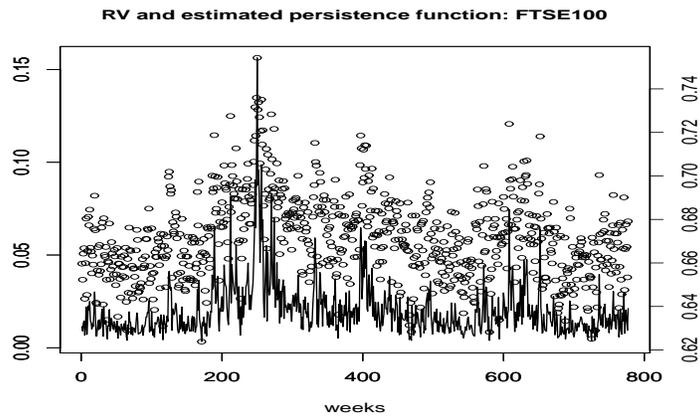}
\caption{Realized volatility (line, left vertical axis) and estimated persistence function (points, right vertical axis). SDAR(M2) for FTSE100.}\label{fig3}
\end{figure}

\subsection{Forecasting accuracy}
Constructing multi-period forecasts for nonlinear models is considerably more difficult than for linear models since exact analytical solutions are not available. For a general discussion about a number of methods of obtaining forecasts for general nonlinear models see Granger and Terasvisrta (1993). Clements and Smith (1997) compare a number of alternative methods of obtaining multi-period forecasts using SETAR models and they conclude that Monte Carlo method performs reasonably well. For this reason, in this paper, SETAR forecasts are generated by Monte Carlo simulation (denoted MC).

In the same way we proceed for the SDAR model. In fact, MC method is a simple simulation method for obtaining multi-step forecasts that can be applied as easily to complex models as the SDAR model. Denote with $H$ the forecast horizon. Let $\tilde{y}^{(m)}_{n+h}$ be the $(n+h)$-th forecast with $h=1,...,H$ obtained in the $m$-th simulation where $n$ is the time of the last observation in the sample. The forecast simulation scheme is then
\begin{equation}
\left\{\begin{array}{cc}
                                       \tilde{y}^{(m)}_{n+1}=\hat{\alpha}+\psi(y_{n};\hat{\boldsymbol{\gamma}})y_{n}+\tilde{\xi}^{(m)}_{n+1} ,\\
                                      \tilde{y}^{(m)}_{n+h}=\hat{\alpha}+\psi(\tilde{y}^{(m)}_{n+h-1};\hat{\boldsymbol{\gamma}})\tilde{y}^{(m)}_{n+h-1}+\tilde{\xi}^{(m)}_{n+h}, &  h=2,...,H\\
                                      \end{array}\right ..
\end{equation}
where $\tilde{\xi}^{(m)}_{n+1}$ and $\tilde{\xi}^{(m)}_{n+h}$ are normally distributed with zero mean and standard deviation $\hat{\sigma}$ for each simulation $m$. Now, averaging these forecasts across the $m=1,..,M$ iterations of the MC yields
$$\tilde{y}_{n+h}=\sum_{m=1}^M \tilde{y}^{(m)}_{n+h}, \quad h=1,..,H$$

To compare the accuracy of the forecasts obtained by SDAR and SETAR models we use three different measures, the Mean Absolute Forecast Error (MAFE), the Mean Square Forecast Error (MSFE) and the Mean Percentage Forecast Error (MAPE). In particular, we compute the Relative Efficiency (RE) measure defined as the ratio of one of the considered measures of the SDAR model and its competitor, the SETAR model. In other words, if we consider the MAFE  $$RE=\frac{MAFE(SDAR)}{MAFE(SETAR)}.$$ A similar computation holds for MSFE and MAPE.
A value of RE lesser or equal than unity indicates that the SDAR model provides more accuracy than the SETAR model.
Results are summarized in tables 5 and figures 4-6.
We consider a forecast horizon from 1 week to 20 weeks ahead (approximately corresponding to a period of 5 months of forecast). Figures 4-6 depict the RE of the three measures of accuracy for each financial index.
We can see that the relative forecast accuracy of the SDAR model with respect to the SETAR(2,4,4) model is steadily improved for the realized volatility extracted from the FTSE100 regardless of the adopted measure of accuracy
(figure 6 for a forecast horizon from 2 to 20 weeks. In particular, the SDAR model is approximately 15-20\% better than the SETAR model for short horizons, i.e., of 10-12 weeks. Differently, in the case of the CAC40 the SDAR model yields superior forecasts up until about 13 weeks ahead after which it gets worse (figure 4). Finally, in the case of the DAX30, the comparison is more complicated since neither model seems to prevail over the other in a systematic way or for defined forecast horizons, even if for short horizons the SETAR(2,2,3) model dominates the SDAR model up to 15\%. On the other hand, for long forecast horizon the RE approaches to 1, i.e., the models equates (figure 5) and often the SDAR model provides more accurate forecasts. Forecasts of 1-week ahead deserve a separate discussion. Here, the SETAR model perform significantly better than the SDAR model for all three indices and independently on the adopted measures. The reason can be found in the fact that for forecasting the $(n+1)$-th observation conditional on $y_n$ the regime is known with certainty, whereas for forecast horizon greater or equal than 2 the regime must be determined by the simulated value $\tilde{y}_{n+1}$ (subjected by an error term).
As regards the measures of accuracy used, MAFE, MSFE and MAPE, the dynamics of the RE in terms of the forecast horizon is very similar. However, we can observe that the MSFE is the one that provides the lowest or highest values of the RE. However, we cannot say that there are appreciable differences between the three measures in terms of a different assessment of the forecasts accuracy.

\begin{table}
\begin{center}
\footnotesize
\begin{tabular}{c|ccc|ccc|ccc}
\hline
   & &  CAC40 & & &  DAX30  &  & &   FTSE100   \\
 \hline
 & MAFE  &  MSFE  &  MAPE &  MAFE  & MSFE  &  MAPE & MAFE  & MSFE  &  MAPE  \\
\hline
$H=1$     &1.5070	&2.2708&	1 .4070&	 1.4117	&  1.9927	&   1.4119	&3.6388&3.9088	&3.6379\\
 $H=2$ &  1.0226  &  0.9861  &  0.9666  &  1.0996 &   1.0914  &  1.0773  &  1.0375  &  0.7343  &  0.9800\\
  $H=3$ &   0.8795   & 0.8217  &  0.8379  &  1.0976 &   1.1212  &  1.0811  &  0.9381  &  0.7546  &  0.9111\\
  $H=4$ &   0.7392  &  0.6631  &  0.7043  &  1.0916 &   1.1661  &  1.0795  &  0.8580  &  0.6808  &  0.8424\\
 $H=5$ &    0.8213  &  0.6925   & 0.7731  &  1.0437 &   1.0100  &  1.0272  &  0.7669  &  0.6258  &  0.7602\\
  $H=6$ &   0.8099  &  0.7217   & 0.7654  &  1.0829 &   1.1312  &  1.0721  &  0.8915  &  0.8119  &  0.8832\\
  $H=7$ &   0.7148  &  0.5848   & 0.6777  &  1.0837 &   1.0704  &  1.0726  &  0.7652  &  0.6096  &  0.7606\\
  $H=8$ &   0.7976  &  0.7754   & 0.7496  &  0.9768 &   0.9448  &  0.9735  &  0.9070  &  0.7813  &  0.8941\\
  $H=9$ &   0.8671  &  0.8612   & 0.8153  &  1.0786 &   1.0061  &  1.0880  &  0.9486  &  0.8191  &  0.9338\\
  $H=10$ &   0.9189  &  0.9072   & 0.8704  &  1.0099 &   1.0074  &  1.0121  &  0.9791  &  0.8740  &  0.9669\\
 $H=11$ &    0.8525  &  0.7846   & 0.8072  &  0.9406 &   0.9542  &  0.9317  &  0.9077  &  0.7769  &  0.8984\\
 $H=12$ &    0.9347  &  0.9234   & 0.8823  &  1.0152 &   0.9615  &  1.0197  &  0.8672  &  0.7368  &  0.8582\\
  $H=13$ &   0.9953  &  0.9806   & 0.9477  &  1.0388 &   1.0396  &  1.0427  &  0.9061  &  0.7600  &  0.8919\\
 $H=14$ &    1.0132  &  1.0852   & 0.9610  &  0.9759 &   0.9911  &  0.9827  &  0.9899  &  0.8675  &  0.9696\\
 $H=15$ &    1.0326  &  1.1621   & 0.9804  &  0.9734 &   0.9373  &  0.9691  &  0.9073  &  0.8076  &  0.8920\\
 $H=16$ &    1.0802  &  1.1701  &  1.0307  &  1.0219 &   0.9886  &  1.0344  &  1.0045  &  0.9968  &  0.9885\\
 $H=17$ &    1.1174  &  1.2889  &  1.0652  &  1.0272 &   1.0171  &  1.0272  &  0.9828  &  0.9080  &  0.9619\\
 $H=18$ &    1.1330  &  1.3335  &  1.0696  &  1.0280 &   1.0127  &  1.0283  &  0.8837  &  0.7906  &  0.8732\\
 $H=19$ &    1.0710  &  1.1556  &  1.0229  &  1.0204 &   1.0083  &  1.0271  &  0.8879  &  0.7600  &  0.8736\\
 $H=20$ &    1.1327 &   1.3331  &  1.0885  &  1.0114 &   0.9881  &  1.0179  &  0.9434  &  0.8996  &  0.9275\\

 \hline

\end{tabular}
\caption{Relative efficiency of the SDAR model in terms of the SETAR model according to the measures of forecast accuracy considered, MAFE, MSFE and MAPE. A value of the ratio lesser than 1 indicates that the SDAR model ensures more accuracy than the SETAR model. We compare SDAR(M1) and SETAR(2,3,3) for CAC40, SDAR(M1) and SETAR(2,2,3) for DAX30, SDAR(M2) and SETAR(2,4,4) for FTSE100.}\label{table5}
\end{center}
\end{table}
\bigskip

\begin{figure}[htbp]
\centering
\includegraphics[width=12cm, height=6cm]{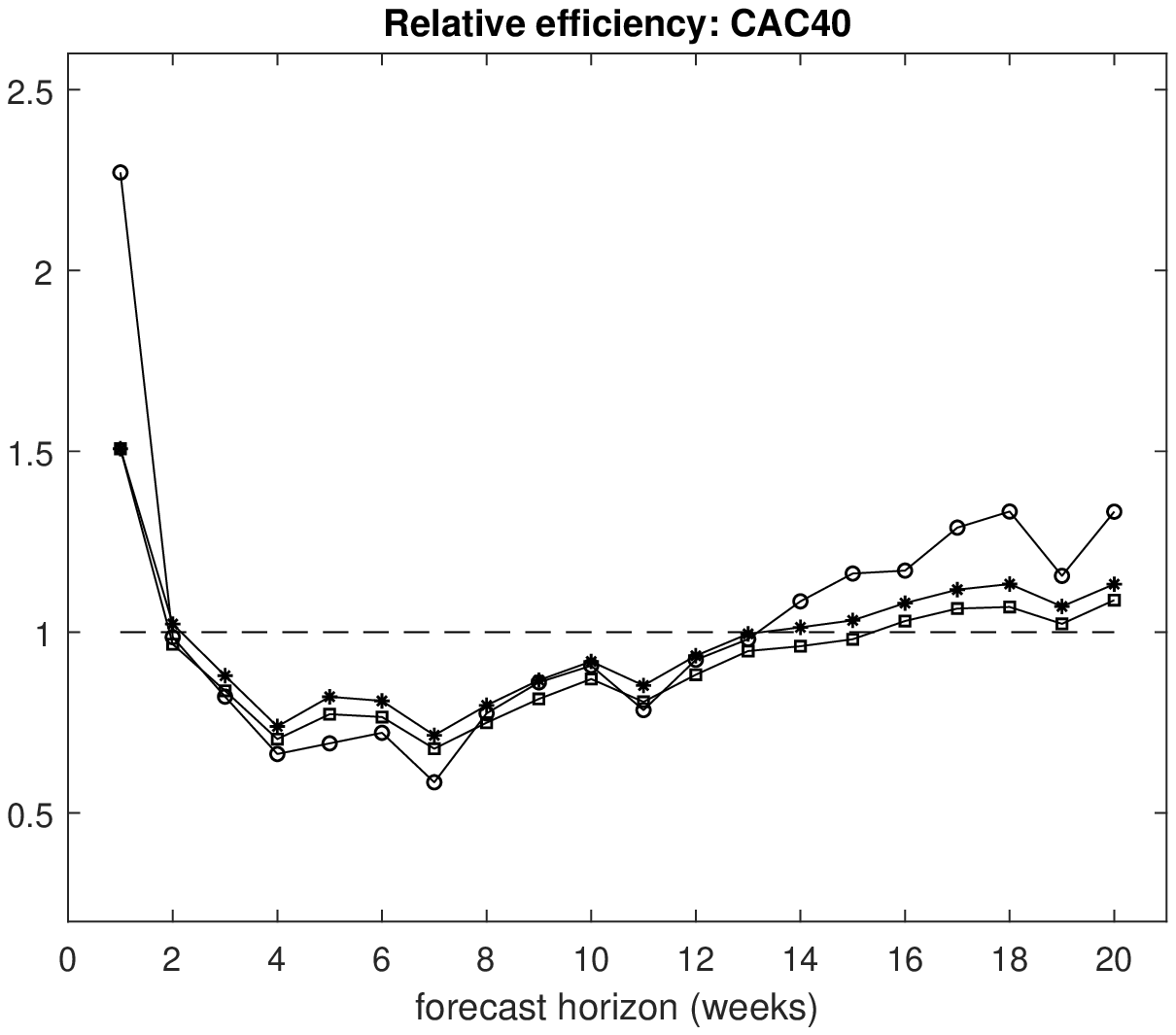}
\caption{Relative efficiency in terms of MAFE (asterisk), MSFE (circle) and MAPE (square). The SDAR model forecast accuracy measure is expressed relative to the corresponding for SETAR(2,3,3) model.}\label{fig4}
\end{figure}

\begin{figure}[htbp]
\centering
\includegraphics[width=12cm, height=6cm]{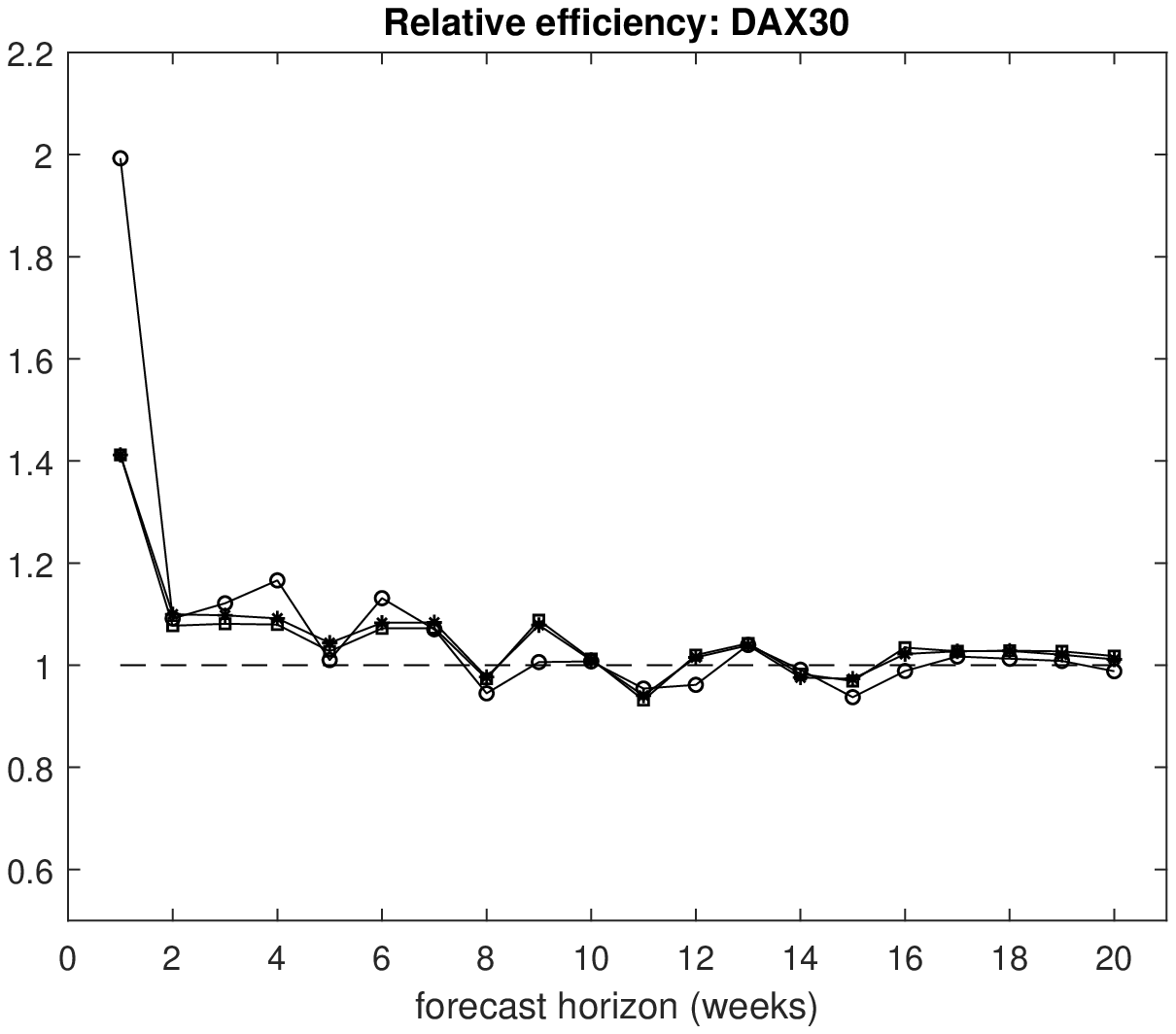}
\caption{Relative efficiency in terms of MAFE (asterisk), MSFE (circle) and MAPE (square). The SDAR model forecast accuracy measure is expressed relative to the corresponding for SETAR(2,2,3) model.}\label{fig5}
\end{figure}

\begin{figure}[htbp]
\centering
\includegraphics[width=12cm, height=6cm]{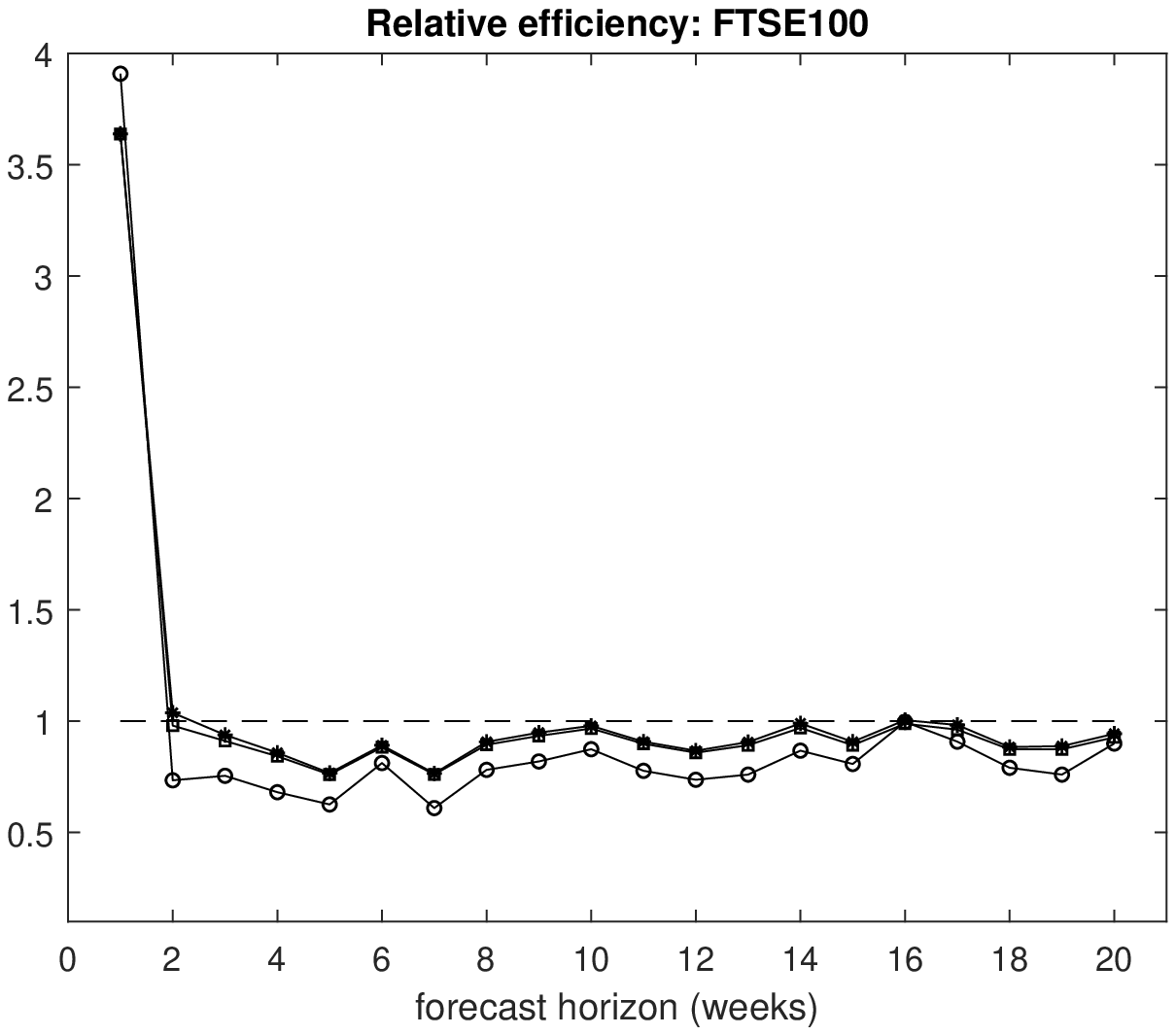}
\caption{Relative efficiency in terms of MAFE (asterisk), MSFE (circle) and MAPE (square).} The SDAR model forecast accuracy measure is expressed relative to the corresponding for SETAR(2,2,4) model.\label{fig6}
\end{figure}

%\begin{figure}[htbp]
%\centering
%\includegraphics[width=12cm, height=6cm]{RE_AR_SDAR_MSE.eps}
%\caption{Relative MSFEs for the SDAR model expressed relative to that for linear AR model. For each index realized %volatility a different AR model is estimated. CAC40 (asterisk), DAX30 (circle) and FTSE100 (square).}\label{fig7}
%\end{figure}

\section{Concluding remarks}
In this paper we provide conditions under which the SDAR model is strictly stationary and uniformly ergodic. The model is a special case of the class of functional-coefficient AR processes in which the autoregressive coefficient is a function of the lagged state variable. We impose a number of assumptions on the persistence function $\psi$ to get nonlinear time series which are strictly stationary and uniformly ergodic. From an estimation point of view, we propose a quasi-maximum likelihood technique and we establish that the estimator is consistent and asymptotically normal. An empirical application to weekly realized volatilities extracted from three European financial indices (CAC40, DAX30 and FTSE100) is presented and the forecast accuracy of the model is discussed relating to an alternative approach among the most used in econometric literature, the SETAR model. We show that the SDAR model has a better predictive ability than the SETAR model in the case of CAC40 and FTSE100 (except for a forecast of 1-step ahead) and just worse in the case of DAX30 for short horizons, offering a possible alternative in modelling and forecasting nonlinear economic time series. It is our belief that further investigations are needed to fully understand the potentiality of SDAR models.

\section{Appendix 1}
In this appendix we shall prove the main results introduced in this paper.
Since in propositions 2.1 e 2.2 ${\boldsymbol\theta}$ is assumed fixed, in order to simplify the notation, in the proofs of these results we  will drop the dependence of $\psi$ on ${\boldsymbol\gamma}$, i.e., we set $\psi(y)=\psi(y;{\boldsymbol\gamma})$.\\

\vspace{0.1cm}
\textbf{Proof of Proposition 2.1}
\vspace{0.1cm}
\begin{proof}
We will apply theorem
4.40 in Douc et al. (2014). Preliminarily, it is
necessary to consider the SDAR model specified in (2.2) as a special case of a more general class of models known as
"iterated random function" satisfying the following recurrence equation $Y_t=f_{\xi_t}(Y_{t-1})$, for all $t \in \mathbb{N}$,
where $f_{\xi_{t}}(Y_{t-1})=\psi(Y_{t-1})Y_{t-1}+\xi_t$.
Theorem 4.40 in Douc et al. (2014) states that, under their assumptions A4.36-A4.38, for
all $y_0\in \mathbb{R}$ the composition $y^{(k)}_t(y_0)=f_{\xi_t} \circ \cdot \cdot \cdot \circ f_{\xi_{t-k-1}}(y_0)$ converges $\mathbb{P}$-a.s. to a r.v. $\tilde{y}_t$ which does not
depend on $y_0$ and, moreover, $(\tilde{y}_t)_t$ is the only strictly stationary solution of $Y_t$.
Therefore, our goal is to verify that under the SDAR model in (2.2) requirements A4.36-A4.38 of Douc et al (2014) are satisfied. We will discuss them separately.
\begin{itemize}
  \item[A4.36]  \emph{The sequence $(\xi_t)_t$ is strict stationary and ergodic}.\\
   Since we assume that the error
term sequence is i.i.d., the requirement trivially holds.

  \item[A4.37] \emph{There exists a measurable function $e\mapsto K_e$ such that $|f_e(y)-f_e(z)|\leq K_e|y-z|$,
$\forall (y, z,e) \in \mathbb{R}^3$, $\mathbb{E}[\ln^+(K_{\xi_0})] < \infty$ and $\mathbb{E}[\ln(K_{\xi_0})] < 0$}.\\
 In our SDAR
model $f_e(y)=\psi(y)y+e$, therefore $|f_e(y)-f_e(z)|=|\psi(y)y-\psi(z)z|$. Assumption \textbf{a1} ensures that $\psi(y)y$ is
a Lipschitz function with Lipschitz constant $L \in (0,1)$ so that $|\psi(y)y-\psi(z)z|\leq
L|y-z|$. Therefore, setting $K_e=L$ for all $e$ all conditions are satisfied.

  \item[A4.38]  \emph{There exists $\hat{y}\in \mathbb{R}$ such that $\mathbb{E}\left[\ln^+(|\hat{y}-f_{\xi_1}(\hat{y})|)\right]<\infty$}.\\
       In the SDAR  model, since $\xi_1\sim N(0,\sigma)$ we can write
      $$\mathbb{E}\left[\ln^+(|\hat{y}-\psi(\hat{y})\hat{y}-\xi_1|)\right]=$$
      $$=\frac{1}{2\sigma \sqrt{2 \pi}}\int_{-\infty}^{+\infty}
      \ln(|\hat{y}-\psi(\hat{y})\hat{y}-s|)\textbf{1}_{\left\{\left\{s<\hat{y}-\psi(\hat{y})\hat{y}-1\right\}\cup \left\{s>\hat{y}-\psi(\hat{y})\hat{y}-1\right\}\right\}} e^{-\frac{s^2}{2 \sigma^2}}ds$$
      which is clearly finite.
\end{itemize}
So, under assumption {\bf a1}, the SDAR model satisfies A4.36-A4.38 in Douc et al (2014), and the result is proved.
\end{proof}

\textbf{Proof of Proposition 2.2}
\vspace{0.1cm}
\begin{proof}
The result is an immediate consequence of proposition 6.8 in Douc et al. (2014) and more specifically of their example 6.11,
taking into account that the density of the kernel of the SDAR model with respect to the Lebesgue measure is given by
$$p(y,\tilde{y})=\frac{1}{\sigma \sqrt{2 \pi}} e^{-\frac{1}{2 \sigma^2} (\tilde{y}-\psi(y)y)^2}.$$
\end{proof}

\textbf{Proof of Theorem 2.1}

\begin{proof} The theorem follows from theorems 3.13 and 6.4 in White (1994) whose assumptions are satisfied by the SDAR model as we are going to show.
Actually, White's requirements 2.1, 2.3, 3.6 and 3.9 are trivially satisfied under our assumptions \textbf{a1}-\textbf{a4}.
So, we analyze all non trivial assumptions 3.1, 3.2$^\prime$, 3.7, 3.8 and 6.1.

\begin{itemize}
\item[3.1]
               \emph{ \textbf{a}) $\mathbb{E}\left[\ell_t(Y^t;\boldsymbol{\theta})\right]$ exists and it is finite for all
                $\boldsymbol{\theta} \in \Theta$;
                \textbf{b}) $\mathbb{E}\left[\ell_t(Y^t;\cdot)\right]$ is continuous on $\Theta$;
                \textbf{c}) $(\ell_t(Y^t;\boldsymbol{\theta}))_t$ satisfies a strong uniform law of large numbers.}
                \\
Notice that $\ell_t(Y^t;\boldsymbol{\theta})=-\ln(\sigma)-\frac{\xi_t^{2}}{2 \sigma^2}=q(\xi_t;\sigma)$. Since  $(\xi^{2}_t)_t$
is an i.i.d. sequence of r.vs. with finite moments requirements \textbf{a}) and \textbf{b}) are satisfied.
In order to verify condition \textbf{c}), we rewrite $q(\xi_t;\sigma)$ as $q(\xi_t;\sigma)=\hat q(z_t;\sigma)=-\ln(\sigma)-\frac{z^2_{t}}{2}$ where
$z_t \sim i.i.d. \ N(0,1)$. Since
$$
\underset{\sigma\in\mathcal S}\sup\left\vert \frac 1n\sum_{t=1}^n\hat q(z_t;\sigma)-\mathbb E\left [\hat q(z_1;\sigma)\right ]\right\vert=
\frac 12\left\vert \frac 1n\sum_{t=1}^nz_t^2-1\right\vert\to 0
$$
condition \textbf{c}) is satisfied.

\item[3.2$^{\prime}$] \emph{ Let $L_n(\boldsymbol{y}^n;\boldsymbol{\theta})=\frac{1}{n}\sum_{t=1}^n \ell_t(Y^t;\boldsymbol{\theta})$ and let $\bar{L}_n(\boldsymbol{\theta})=\mathbb{E}[L_n(\boldsymbol{y}^n;\boldsymbol{\theta})]$. Assume that
$\bar{L}_n(\boldsymbol{\theta})$ is $O(1)$ uniformly in $\Theta$}.\\
It is an immediate consequence of assumption \textbf{a3}.

\item[3.7]\emph{\textbf{a}) $\mathbb{E}\left[\frac{1}{n} \sum_{t=1}^n \nabla_{\boldsymbol{\theta}} \ell_t(Y^t;\boldsymbol{\theta})\right]<+\infty$ for all
                $\boldsymbol{\theta} \in \Theta$, uniformly in $n$;
               \textbf{b}) $\mathbb{E}\left[\frac{1}{n} \sum_{t=1}^n \nabla_{\boldsymbol{\theta}}  \ell_t(Y^t;\boldsymbol{\theta})\right]$ is continuous on $\Theta$ uniformly in $n$;
                \textbf{c}) the sequence $(\nabla_{\boldsymbol{\theta}}  \ell_t(Y^t;\boldsymbol{\theta}))_t$ satisfies a strong uniform law of large numbers}.\\
 Notice that the elements of the gradient $\nabla_{\boldsymbol{\theta}}  \ell_t(Y^t;\boldsymbol{\theta})$ are two types of sequences of r.vs.: $\ell_t^{\alpha}(Y^t;\boldsymbol{\theta})=\frac{\xi_t}{\sigma^2}$ and $\ell_t^{\sigma}(Y^t;\boldsymbol{\theta})=\frac{1}{\sigma^3}(\xi_t^{2}-\sigma^2)$ are both i.i.d. sequences, whereas $\ell_t^{\gamma_k}(Y^t;\boldsymbol{\theta})=\frac{1}{\sigma^2}\psi^{\gamma_k}(Y_{t-1};\boldsymbol{\gamma})\xi_t Y_{t-1}$ is a martingale difference sequence, for all $k=1,...,p$.
 Requirements \textbf{a}) and \textbf{b}) follow directly from our assumptions. As for condition \textbf{c}), let us start considering
 the sequence
 $\ell_t^{\alpha}(Y^t;\boldsymbol{\theta})$. Following the same idea used in the proof of condition \textbf{c}) of assumption 3.1 in White (1994) above,
 we have
 $$\underset{\boldsymbol{\theta}\in\Theta}\sup\left\vert\frac 1n\sum_{t=1}^n
\ell_t^{\alpha}(Y^t;\boldsymbol{\theta})-\mathbb E[\ell_t^{\alpha}(Y^t;\boldsymbol{\theta})]\right\vert=$$$$=
 \underset{\sigma\in\mathcal S}\sup\frac 1\sigma\left\vert \frac 1n\sum_{t=1}^nz_t\right\vert=
\frac 1{\sigma_{min}}\left\vert \frac 1n\sum_{t=1}^nz_t\right\vert\to 0
$$
where $\sigma_{min}=\min\mathcal S>0$ by assumption \textbf{a3}.
 As for $\ell_t^{\sigma}(Y^t;\boldsymbol{\theta})$ we can use the same arguments. Regarding the sequence $(\ell_t^{\gamma_k}(Y^t;\boldsymbol{\theta}))_t$ which has a more complex serial dependence structure, we need to exploit theorem 1 in Potscher and Prucha (1989) where the authors state a strong uniform law of large numbers for functionals relative to dependent and heterogeneous sequence of r.vs.. We need to control assumptions 1, 2, 3, 4 and 5A in the mentioned paper.
Most assumptions are trivially satisfied under the hypotheses of the SDAR model.
The only non trivial ones are assumption 4, for which $(\ell_t^{\gamma_k}(Y^t;\boldsymbol{\theta}))_t$ is required to satisfy
a pointwise strong law of large numbers and assumption 3 which states that, for $k=1,\ldots ,p$,
\begin{equation}\label{dominance}\sup_n \frac{1}{n} \sum_{t=1}^n
\mathbb{E}\left[\sup_{\boldsymbol{\theta}}\left|\ell_t^{\gamma_k}(Y^t;\boldsymbol{\theta}) \right|^{1+\delta} \right]< \infty.\end{equation}
As for the first one, thanks to assumption {\bf a4} and Remark 2.2, the conclusion immediately follows from Theorem 3.77 in White (1984).
Moreover, (5.4) is verified since, under assumption {\bf a5},
$$\begin{aligned} \mathbb{E}\left[\sup_{\boldsymbol{\theta}}\left|\ell_t^{\gamma_k}(Y^t;\boldsymbol{\theta}) \right|^{1+\delta}\right]
& =\mathbb{E}\left[\sup_{\boldsymbol{\theta}}\left|\frac{1}{\sigma^2}\psi^{\gamma_k}(Y_{t-1};\boldsymbol{\gamma})\xi_t Y_{t-1} \right|^{1+\delta}\right]\leq \\
& \leq \frac{C^{1+\delta}}{\sigma_{\min}^{1+\delta}} \mathbb{E}\left[\vert Z\vert^{1+\delta}\right] < \infty,
\end{aligned}$$
where $Z\sim N(0,1)$.

\item[3.8]               \emph{\textbf{a}) $\mathbb{E}\left[\frac{1}{n} \sum_{t=1}^n \nabla_{\boldsymbol{\theta}} ^{2} \ell_t(Y^t;\boldsymbol{\theta})\right]<+\infty$ for all
                $\boldsymbol{\theta} \in \Theta$, uniformly in $n$;
                 \textbf{b}) $\mathbb{E}\left[\frac{1}{n} \sum_{t=1}^n \nabla_{\boldsymbol{\theta}}^{2} \ell_t(Y^t;\boldsymbol{\theta})\right]$ is continuous on $\Theta$ uniformly in $n$;
                \textbf{c}) the sequence $(\nabla_{\boldsymbol{\theta}}^{2} \ell_t(Y^t;\boldsymbol{\theta}))_t$ satisfies a strong uniform law of large numbers}.\\
The elements of the hessian matrix $\nabla_{\boldsymbol{\theta}}^{2} \ell_t(Y^t;\boldsymbol{\theta})$ are the following: $$\ell_t^{\alpha \alpha}(Y^t;\boldsymbol{\theta})=-\frac{1}{\sigma^2},$$
$$\ell_t^{\alpha \gamma_k}(Y^t;\boldsymbol{\theta})=-Y_{t-1}\frac{\psi^{\gamma_k}(Y_{t-1};\boldsymbol{\theta})}{\sigma^2},
\quad for \ all \ k,$$
$$\ell_t^{\alpha \sigma}(Y^t;\boldsymbol{\theta})=-\frac{2 \xi_t}{\sigma^3},$$
$$\ell_t^{\gamma_k \gamma_k}(Y^t;\boldsymbol{\theta})=\xi_t Y_{t-1}\frac{\psi^{\gamma_k \gamma_k}(Y_{t-1};\boldsymbol{\theta})}{\sigma^2}-Y^2_{t-1}\frac{(\psi^{\gamma_k}(Y_{t-1};
\boldsymbol{\gamma}))^2}{\sigma^2},\quad for \ all \ k,$$
$$\ell_t^{\gamma_k \gamma_j}(Y^t;\boldsymbol{\theta})=$$$$=\xi_t Y_{t-1}\frac{\psi^{\gamma_k \gamma_j}(Y_{t-1};\boldsymbol{\gamma})}{\sigma^2}-Y^2_{t-1}
\frac{\psi^{\gamma_k}(Y_{t-1};\boldsymbol{\gamma})\psi^{\gamma_j}(Y_{t-1};\boldsymbol{\gamma})}{\sigma^2},
 \quad for \ all \ k \geq j,$$
$$\ell_t^{\gamma_k \sigma}(Y^t;\boldsymbol{\theta})=-\xi_t Y_{t-1} \frac{2\psi^{\gamma_k}(Y_{t-1};\boldsymbol{\gamma})}{\sigma^3}, \quad for \ all \ k, $$
and
$$\ell_t^{\sigma \sigma}(Y^t;\boldsymbol{\theta})=\frac{1}{\sigma^2}-\frac{3\xi_t^{2}}{\sigma^4}.$$
Notice that we can make use of similar arguments as those used in the previous points to establish that each of these sequences satisfy a strong uniform law of large numbers. Indeed, we are dealing with i.i.d. sequences ($(\ell_t^{\alpha \alpha}(Y^t;\boldsymbol{\theta}))_t$ and $(\ell_t^{\sigma \sigma}(Y^t;\boldsymbol{\theta}))_t$), stationary and ergodic sequences ($(\ell_t^{\alpha \gamma_k}(Y^t;\boldsymbol{\theta}))_t$, for all $k$) and martingale difference sequences ($(\ell_t^{\gamma_k \sigma}(Y^t;\boldsymbol{\theta}))_t$, for all $k$), whereas
$(\ell_t^{\gamma_k \gamma_k}(Y^t;\boldsymbol{\theta}))_t$ and $(\ell_t^{\gamma_k \gamma_j}(Y^t;\boldsymbol{\theta}))_t$ are linear combinations of them.  We will only discuss the last case by considering the others
analogous.
Thanks to assumption \textbf{a4} and Remark 2.2, Theorems 3.34 and 3.77 in White (1984) ensure that a pointwise strong law of large numbers is satisfied by $\ell_t^{\gamma_k \gamma_k}(Y^t;\boldsymbol{\theta})$ and by $\ell_t^{\gamma_k \gamma_j}(Y^t;\boldsymbol{\theta})$.
The most delicate issue is to check condition 3 in Potscher and Prucha (1989). We consider the sequence $\ell_t^{\gamma_k\gamma_j}(Y^t;\boldsymbol{\theta})$ being
$\ell_t^{\gamma_k\gamma_k}(Y^t;\boldsymbol{\theta})$ completely equivalent. We remark that, under assumption \textbf{a5},
$$\begin{aligned}\mathbb{E}&\left[\sup_{\boldsymbol{\theta}}\left|\ell_t^{\gamma_k \gamma_j}(Y^t;\boldsymbol{\theta}) \right|^{1+\delta}\right]=\\
&=\mathbb E\left [\sup_{\boldsymbol{\theta}}\left \vert \xi_t Y_{t-1}\frac{\psi^{\gamma_k \gamma_j}(Y_{t-1};\boldsymbol{\gamma})}{\sigma^2}-Y^2_{t-1}\frac{\psi^{\gamma_k}(Y_{t-1};\boldsymbol{\gamma})\psi^{\gamma_j}(Y_{t-1};\boldsymbol{\gamma})}{\sigma^2}\right\vert^{1+\delta}\right ]\leq\\
&\leq \left [\frac{ D}{\sigma_{min}}\left (\mathbb E[\vert Z\vert^{1+\delta}]\right )^{\frac 1{1+\delta}}+
\frac{ C^2}{\sigma_{min}^2}\right ]^{1+\delta} < \infty
  \end{aligned}
$$
by Minkowski's inequality.

\item[6.1] \emph{The array $\frac{1}{\sqrt{n}} s^0_{t}=\frac{1}{\sqrt{n}}\nabla_{\boldsymbol{\theta}}  \ell_t(Y^t;\boldsymbol{\theta}^0)$ obeys the central limit theorem with covariance matrix
$$B^0_{n}=var(\frac{1}{\sqrt{n}}\sum_{t=1}^n s^0_{t})=\mathbb{E}[\frac{1}{n}
\sum_{t=1}^n \nabla_{\boldsymbol{\theta}}  \ell_t(Y^t;\boldsymbol{\theta}^0)\left(\nabla_{\boldsymbol{\theta}}  \ell_t(Y^t;\boldsymbol{\theta}^0)\right)^T ],$$ which is $O(1)$ and uniformly positive definite}.\\
We prove this statement relative to the SDAR model thanks to Crowder (1976). We refer
to conditions 3.6 and 3.7 of that paper. Set $X_t= \mathbf{1}^{T} \nabla_{\boldsymbol{\theta}}  \ell_t(Y^t;\boldsymbol{\theta}^0)$.
Notice that $(X_t)_t$ is a martingale difference sequence, i.e., $\mathbb{E}_{t-1}\left[X_t\right]=0$, $\mathbb{P}^0-$a.s., where $\mathbb{E}_{t-1}$ denotes the conditional expectation $\mathbb{E}[\cdot| Y_{t-1}]$.
Condition 3.6 in Crowder (1976) is
$$\frac{\sum_{t=1}^n \mathbb{E}_{t-1} \left[X_t^{2}\right]}{\sum_{t=1}^n \mathbb{E} \left[X_t^{2}\right]} \stackrel{\mathbb{P}}{\longrightarrow}1.$$
But
$$\mathbb{E}_{t-1} r\left[X_t^{2}\right]=\mathbb{E}_{t-1} \left[(\mathbf{1}^{T} \nabla_{\boldsymbol{\theta}} \ell_t(Y^t;\boldsymbol{\theta}^0))^2\right]=$$
$$=\mathbb{E}_{t-1}\left[\left(\ell_t^{\alpha}(Y^t;\boldsymbol{\theta}^0)+\sum_{k=1}^p \ell_t^{\gamma_k}(Y^t;\boldsymbol{\theta}^0)+
\ell_t^{\sigma}(Y^t;\boldsymbol{\theta}^0)\right)^2\right]=$$
$$
=\mathbb{E}_{t-1}\left[\left(\ell_t^{\alpha}(Y^t;\boldsymbol{\theta}^0)
\right)^2\right]+\mathbb{E}_{t-1}\left[\left(\sum_{k=1}^p \ell_t^{\gamma_k}(Y^t;\boldsymbol{\theta}^0)
\right)^2\right]+$$
$$
+\mathbb{E}_{t-1}\left[\left(\ell_t^{\sigma}(Y^t;\boldsymbol{\theta}^0)
\right)^2\right]+2\mathbb{E}_{t-1}\sum_{k=1}^p \ell_t^{\alpha}(Y^t;\boldsymbol{\theta}^0)\ell_t^{\gamma_k}(Y^t;\boldsymbol{\theta}^0)+$$$$+2\mathbb{E}_{t-1}\sum_{k\neq j} \ell_t^{\gamma_k}(Y^t;\boldsymbol{\theta}^0)
\ell_t^{\gamma_j}(Y^t;\boldsymbol{\theta}^0)+2\mathbb{E}_{t-1}\ell_t^{\alpha}(Y^t;\boldsymbol{\theta}^0) \ell_t^{\sigma}(Y^t;\boldsymbol{\theta}^0),
$$
since it is easy to verify that $\mathbb{E}_{t-1}\left[\ell_t^{\sigma}(Y^t;\boldsymbol{\theta}^0) \ell_t^{\gamma_k}(Y^t;\boldsymbol{\theta}^0)\right]=0$, for all $k=1,..,p$, and
$\mathbb{E}_{t-1}\left[\ell_t^{\alpha}(Y^t;\boldsymbol{\theta}^0) \ell_t^{\sigma}(Y^t;\boldsymbol{\theta}^0)\right]=0$.
Substituting the expressions of partial derivatives of $\ell_t(Y^t;\boldsymbol{\theta})$ we get

$$\mathbb{E}_{t-1} \left[X_t^{2}\right]=\frac{3}{\sigma^2}+\frac{2 Y_{t-1}}{\sigma^2}\sum_{k=1}^p \psi^{\gamma_k}(Y_{t-1};\boldsymbol{\gamma}^0)+$$$$+\frac{Y^2_{t-1}}{\sigma^2}\left[
\sum_{k =1 }^{p} \left(\psi^{\gamma_k}(Y_{t-1};\boldsymbol{\gamma}^0)\right)^2+2\sum_{k \neq j}
\psi^{\gamma_k}(Y_{t-1};\boldsymbol{\gamma}^0)\psi^{\gamma_j}(Y_{t-1};\boldsymbol{\gamma}^0)\right].
$$
We conclude that $\mathbb{E}_{t-1}\left[X_t^{2}\right]$ is strictly stationary and uniformly ergodic sequence of r.vs.. By assumption \textbf{a4} and Remark 2.2 we can apply theorem 3.34 in White (1984) to get
$$\frac{\sum_{t=1}^n \mathbb{E}_{t-1} \left[X_t^{2}\right]}{\sum_{t=1}^n \mathbb{E} \left[X_t^{2}\right]} =
\frac{\frac{1}{n}\sum_{t=1}^n \mathbb{E}_{t-1} \left[X_t^{2}\right]-\frac{1}{n}\sum_{t=1}^n \mathbb{E} \left[X_t^{2}\right]}{\frac{1}{n}\sum_{t=1}^n \mathbb{E} \left[X_t^{2}\right]}+1
\stackrel{\mathbb{P}}{\longrightarrow} 1,$$
since $\mathbb{E}\left\{\mathbb{E}_{t-1}\left[X^2_{t}\right]\right\}=\mathbb{E}\left[X^2_{t}\right]$.\\
Condition 3.7 in Crowder (1976) requires that for some $\delta >0$ and for every $\epsilon >0$,  $\frac{\sum_{t=1}^n\mathbb{E}\left[\left|X_t\right|^{2+\delta}\textbf{1}_{|X_t|\geq \epsilon v_n}\right]}{v_n^{2+\delta}} \rightarrow 0$ where $v_n^{2}=\sum_{t=1}^{n}\mathbb{E}\left[X_t^{2}\right]$.
Choosing $\delta=2$ and considering the stationarity of $\mathbb{E}_{t-1}\left[X_t^{2}\right]$
we get
$$\frac{\sum_{t=1}^n \mathbb{E}\left [X_t^{4} \mathbf{1}_{|X_t|> \epsilon v_n}\right ]}{v_n^{4}} \leq
\frac{\sum_{t=1}^n \mathbb{E}\left [X_t^{4}\right ]}{v_n^{4}} \leq \eta\cdot \frac 1n \longrightarrow 0$$
where $0<\eta<\infty$ is the upper bound of $\mathbb{E}\left [X_t^{4}\right]$.

Since all required assumptions are satisfied, the result follows from Theorems 3.13 and 6.4 in White (1994).
\end{itemize}
\end{proof}

\section{Appendix 2}
In this appendix we give a brief description of SETAR($p, d_1,d_2$) models used as competitors for evaluating forecasting accuracy of SDAR models. For a detailed discussion on SETAR models the reader can consult Tong (1978, 1980, 1986 and 1995). SETAR models assume that a variable $y_t$ is a linear AR within a regime but may move among regimes depending on the value assumed by the lagged variable $y_{t-d}$. In our application $d=1$ the number of regimes $p$ is 2. More in particular, the process $y_t$ follows an AR($d_1$) process when it is in the "low regime" and an AR($d_2$) process when it is in the "high regime". Formally,
\begin{equation}
\left\{\begin{array}{cc}
                                       y_t=c_{1}+\sum_{i=1}^{d_1} \phi_{1,i} y_{t-i}+\epsilon_{1,t}, \quad y_{t-1} \leq r\\
                                      y_t=c_{2}+\sum_{i=1}^{d_2} \phi_{2,i} y_{t-i}+\epsilon_{2,t}, \quad y_{t-1} > r,\\
                                      \end{array}\right ..
\end{equation}
 where $\epsilon_{j,t} \sim N(0,\sigma_j)$, $j=1,2$.
We report the output of R (package "tsDyn") relating to the estimated SETAR models for the three different index log-volatilities.
\footnotesize
\begin{itemize}
\item CAC40. R output: Non linear autoregressive model\\
\vspace{0.1cm}
SETAR model ( 2 regimes)\\
\vspace{0.1cm}
Coefficients:\\
\vspace{0.1cm}
\begin{tabular}{cccccc}
\hline
Low regime: & & & \\
const.L  &   phiL.1  &   phiL.2 &    phiL.3\\
-1.3292280 & 0.2297850&  0.2442498 & 0.1805137\\
\hline
High regime: & & & \\
   const.H  &   phiH.1  &   phiH.2  &   phiH.3\\
-0.6589916  & 0.4226148 & 0.3088542  & 0.1112698\\
\hline
\end{tabular}\\
\vspace{0.2cm}
Threshold Variable: Z(t) = + (0) X(t)+ (1)X(t-1)+ (0)X(t-2)\\
\vspace{0.2cm}
Proportion of points in low regime: 35.66\%, High regime: 64.34\%\\
\vspace{0.2cm}
\begin{tabular}{cccccc}
\hline
Residuals:   &  &  &\\
 Min     &    1Q   &     Median    &     3Q  &      Max\\
 -1.9040632 & -0.2957076 & -0.0064422 & 0.2994092 &  1.4948586\\
\hline
\end{tabular}\\
\vspace{0.2cm}
Fit: residuals variance = 0.1983,  AIC = -1239, MAPE = 9.375\%\\
\vspace{0.4cm}
\begin{tabular}{ccccccc}
\hline
Coefficient(s)  &  &  & & \\
  &       Estimate & Std. Error & t value  &$Pr(>|t|)$\\
const.L &-1.329228 &   0.454556 & -2.9242& 0.0035545 **\\
phiL.1 &  0.229785 &   0.061935 &  3.7101& 0.0002221 ***\\
phiL.2 &  0.244250  &  0.089647  & 2.7246& 0.0065847 **\\
phiL.3 &  0.180514   & 0.058308 &  3.0959& 0.0020337 **\\
const.H& -0.658992   & 0.189383 & -3.4797& 0.0005303 ***\\
phiH.1 &  0.422615   & 0.043588 &  9.6956& $<$ 2.2e-16 ***\\
phiH.2 &  0.308854   & 0.062391 &  4.9503& 9.112e-07 ***\\
phiH.3 &  0.111270   & 0.045072 &  2.4687 &0.0137753 *\\
\hline
\end{tabular}\\
\vspace{0.2cm}
Signif. codes:  0 ‘***’ 0.001 ‘**’ 0.01 ‘*’ 0.05 ‘.’ 0.1 ‘ ’ 1\\
\vspace{0.2cm}

\begin{tabular}{cc}
  \hline
 BDS Test ($p$-values)\\
  \hline
    eps[1] m=2: 0.5944\\
    eps[1] m=3: 0.6796\\
    eps[2] m=2: 0.8801\\
    eps[2] m=3: 0.9139\\
    eps[3] m=2: 0.9262\\
    eps[3] m=3: 0.9241\\
    eps[4] m=2: 0.9720\\
    eps[4] m=3: 0.9621\\

 \hline

\end{tabular}\\

\item DAX30. R output: Non linear autoregressive model\\
\vspace{0.1cm}
SETAR model ( 2 regimes)\\
\vspace{0.1cm}
Coefficients:\\
\vspace{0.1cm}
\begin{tabular}{cccccc}
\hline
Low regime: & & & \\
const.L  &   phiL.1  &   phiL.2 \\
-1.8719247 &  0.2265772 & 0.2711836\\
\hline
High regime: & & & \\
   const.H  &   phiH.1  &   phiH.2  &   phiH.3\\
-0.8432665 & 0.3662158 & 0.2818168 & 0.1424999 \\
\hline
\end{tabular}\\
\vspace{0.2cm}
Threshold Variable: Z(t) = + (0) X(t)+ (1)X(t-1)+ (0)X(t-2)\\
\vspace{0.2cm}
Proportion of points in low regime: 15.12\%, High regime: 84.88\%\\
\vspace{0.2cm}
\begin{tabular}{cccccc}
\hline
Residuals:   &  &  &\\
 Min     &    1Q   &     Median    &     3Q  &      Max\\
 -1.953924 & -0.260304 & 0.035248  & 0.305523  &1.660976\\
\hline
\end{tabular}\\
\vspace{0.2cm}
Fit: residuals variance = 0.2086,  AIC = -1202, MAPE = 9.367\%\\
\vspace{0.4cm}
\begin{tabular}{ccccccc}
\hline
Coefficient(s)  &  &  & & \\
    &     Estimate & Std. Error & t value  &$Pr(>|t|)$\\
const.L &-1.871925  &  0.679896&  -2.7533 &0.0060397 **\\
phiL.1  & 0.226577  &  0.083477&   2.7142 &0.0067914 **\\
phiL.2  & 0.271184   & 0.142042&   1.9092 &0.0566112 .\\
const.H &-0.843266   & 0.168321&  -5.0099 &6.760e-07 ***\\
phiH.1  & 0.366216   & 0.039314&   9.3152 &< 2.2e-16 ***\\
phiH.2  & 0.281817   & 0.051639 &  5.4574 &6.517e-08 ***\\
phiH.3  & 0.142500   & 0.040288  & 3.5371 &0.0004288 ***\\
\hline
\end{tabular}\\
\vspace{0.2cm}
Signif. codes:  0 ‘***’ 0.001 ‘**’ 0.01 ‘*’ 0.05 ‘.’ 0.1 ‘ ’ 1\\
\vspace{0.2cm}

\begin{tabular}{cc}
  \hline
 BDS Test ($p$-values)\\
  \hline
   eps[1] m=2: 0.8822\\
    eps[1] m=3: 0.7256\\
    eps[2] m=2: 0.5463\\
    eps[2] m=3: 0.4947\\
    eps[3] m=2: 0.2789\\
    eps[3] m=3: 0.2962\\
    eps[4] m=2: 0.1310\\
    eps[4] m=3: 0.1758 \\

 \hline

\end{tabular}\\

\item FTSE100. R output: Non linear autoregressive model\\
\vspace{0.1cm}
SETAR model ( 2 regimes)\\
\vspace{0.1cm}
Coefficients:\\
\vspace{0.1cm}
\begin{tabular}{cccccc}
\hline
Low regime: & & & \\
const.L  &   phiL.1  &   phiL.2 &  \\
-2.1761709 & 0.2340485 & 0.2454682 \\
\hline
High regime: & & & \\
   const.H  &   phiH.1  &   phiH.2  &   phiH.3 & phiH.4\\
-0.5842767 & 0.4288785 & 0.2148472 & 0.0202444 & 0.2013528\\
\hline
\end{tabular}\\
\vspace{0.2cm}
Threshold Variable: Z(t) = + (0) X(t)+ (1)X(t-1)+ (0)X(t-2)\\
\vspace{0.2cm}
Proportion of points in low regime: 35.06\%, High regime: 64.94\%\\
\vspace{0.2cm}
\begin{tabular}{cccccc}
\hline
Residuals:   &  &  &\\
 Min     &    1Q   &     Median    &     3Q  &      Max\\
-1.388613 &-0.298019 & 0.041845 & 0.302668 & 1.377064\\
\hline
\end{tabular}\\
\vspace{0.2cm}
Fit: residuals variance = 0.1968,  AIC = -1245, MAPE = 8.932\%\\
\vspace{0.4cm}
\begin{tabular}{ccccccc}
\hline
Coefficient(s)  &  &  & & \\
    &     Estimate & Std. Error & t value  &$Pr(>|t|)$\\
const.L& -2.176171  &  0.491236&  -4.4300 &1.079e-05 ***\\
phiL.1 &  0.234048  &  0.061853 &  3.7840 &0.0001664 ***\\
phiL.2 &  0.245468  &  0.095003 &  2.5838 &0.0099554 **\\
const.H& -0.584277  &  0.190443 & -3.0680& 0.0022306 **\\
phiH.1 &  0.428879 &   0.043040 &  9.9647 &< 2.2e-16 ***\\
phiH.2 &  0.214847  &  0.064563 &  3.3277 & 0.0009173 ***\\
phiH.3 &  0.020244  &  0.046275 &  0.4375 &0.6618838\\
phiH.4 &  0.201353  &  0.044348 &  4.5403 &6.518e-06 ***\\

\hline
\end{tabular}\\
\vspace{0.2cm}
Signif. codes:  0 ‘***’ 0.001 ‘**’ 0.01 ‘*’ 0.05 ‘.’ 0.1 ‘ ’ 1\\
\vspace{0.2cm}

\begin{tabular}{cc}
  \hline
 BDS Test ($p$-values)\\
  \hline
   eps[1] m=2: 0.5739\\
    eps[1] m=3: 0.8881\\
    eps[2] m=2: 0.8580\\
    eps[2] m=3: 0.6026\\
    eps[3] m=2: 0.8120\\
    eps[3] m=3: 0.7970\\
    eps[4] m=2: 0.6593\\
    eps[4] m=3: 0.8517\\

 \hline

\end{tabular}\\
\end{itemize}


\begin{thebibliography}{}
\bibitem{} Andrews D. K .W. (1988). Laws of Large Numbers for Dependent Non-Identically Distributed Random Variables, \emph{Econometric Theory}, \emph{4(3)}, 458-467.

\bibitem{} Bierens H.J. (1981). Robust Methods and Asymptotic Theory in Nonlinear Econometrics, \emph{Lecture Notes in Economics and Mathematical Systems}, 192. Berlin: Springer Verlag.

\bibitem{} Bierens H.J. (1984). Model Specification Testing of Time Series Regression", \emph{Journal of Econometrics}, \emph{26}, 323-353.

\bibitem{} Brock W.A., Dechert W.D., Sheinkman J.A. (1987). A Test of Independence Based on the Correlation Dimension, SSRI no. 8702, Department of Economics, University of Wisconsin, Madison.

\bibitem{} Brock W.A., Dechert W.D., Sheinkman J.A., LeBaron B. (1996). A Test of Independence Based on the Correlation Dimension, \emph{Econometric Reviews}, \emph{15}, 197-235.

\bibitem{} Cai Z., Fan J., Yao Q. (2000). Functional-Coefficient Regressive Models for Nonlinear Time Series, \emph{J. Am. Statist. Assoc.}, \emph{95(451)}, 941-956.

    \bibitem{} Chappel D., Padmore J., Mistry P., Ellis C. (1996). A threshold model for the Frenchfranc/Deutschmark exchange rate, \emph{Journal of Forecasting}, \emph{15}, 155–164.

%\bibitem{}   Chan, K., Tsay, R. S.(1998): “Limiting Properties of the Least Squares Estimator of a Continuous Threshold Autoregressive Model”, \emph{Biometrika}, \emph{85}, 413-426.

\bibitem{} Chan K.S. (1990). Percentage points of likelihood ratio tests for threshold autoregression, \emph{Journal of Royal Statistical Society B},\emph{53(3)}, 691-696.

     \bibitem{} Chen R., Liu M.L. (2001). Functional coefficient autoregressive models: estimation and tests of hypotheses, \emph{Journal of Time Series Analysis}, \emph{22(2)}, 151-173.

\bibitem{} Chen R., Tsay R. (1993). Functional-coefficient autoregressive models, \emph{J. Am. Statist. Assoc.}, \emph{88}, 298-308.

\bibitem{} Cherubini U., Gobbi F. (2013). A Convolution-based
Autoregressive Process, in F. Durante, W.
 Haerdle, P. Jaworski editors. Workshop on Copula
 in Mathematics and Quantitative Finance. Lecture Notes in
 Statistics-Proceedings. Springer, Berlin/Heidelberg.

\bibitem{} Cherubini U., Gobbi F., Mulinacci S. (2016).
\emph{Convolution Copula Econometrics}, SpringerBriefs in
Statistics.

\bibitem{} Clements M.P., Krozling H.M. (1998). A comparison of the forecast performance of Markov-Swtching and Threshold autoregressive models of US GNP, \emph{The Econometrics Journal}, \emph{1(1)}, C47–C75.

\bibitem{} Clements M.P., Smith J. (1997). The performance of alternative forecasting methods for SETAR
models, \emph{International Journal of Forecasting}, \emph{13}, 463–475.

\bibitem{} Crowder M. J. (1976). Maximum likelihood estimation for dependent observations, \emph{J. Roy. Statist. Soc. Ser.}, \emph{B 38}, 45-53.

%\bibitem{} De Gooijer J.G., Kumar K. (1992): "Some recent developments in non-linear time series modelling, testing %and forecasting", \emph{International journal of Forecasting},
%\emph{8}, 135–156.

%\bibitem{} De Gooijer J.G., De Bruin P. (1997): "On SETAR forecasting", \emph{Statistics and Probability Letters},
%\emph{37}, 7–14.

\bibitem{} Diaconis P., Freedman D.(1999). Iterated Random Functions, \emph{SIAM Rev.}, \emph{47(1)}, 45-76.

%\bibitem{} Diebold F.X., Nason J.A.(1990): "Nonparametric exchange rate prediction", \emph{Journal of International %Economics}, 28, 315-332.

\bibitem{} Douc R., Moulines E., Stoffer D. S. (2014). \emph{Nonlinear Time Series. Theory, Methods and Applications with R Examples}, CRC Press.

\bibitem{} Fan J., Yao Q. (2003). \emph{Nonlinear Time Series: Nonparametric and Patrametric Methods}, Springer-Verlag, New York.

\bibitem{} Gobbi F., Mulinacci S. (2019). Mixing and moments properties of a non-stationary copula-based Markov process, forthcoming in \emph{Communications in Statistics – Theory and Methods}, available online: https://www.tandfonline.com/doi/full/10.1080/03610926.2019.1602653

\bibitem{} Granger C.W.J., Terasvirta T. (1993). \emph{Modelling Nonlinear Economics Relationships}, Oxford University Press, Oxford.

\bibitem{} Haggan V., Ozaki T. (1981). Modelling nonlinear random vibrations using an amplitude-dependent autoregressive time series model,
\emph{Biometrica}, \emph{68(1)}, 189-196.

\bibitem{} Hardle W., Lutkepohl H., Chen R. (1997). A review of nonparametric time series analysis", \emph{International Statistical Review}, 65(1), 49-72.

%\bibitem{} Hansen B. E. (1997): "Inference in TAR Models", \emph{Studies in Nonlinear Dynamics and Econometrics}, \emph{2(1)}, The MIT Press.

\bibitem{} Hastie T., Tibshirani R. (1990). \emph{Generalized additive models}. Chapman \& Hall, New York.


%\bibitem{}  Keenan D. M. (1985): "A Tukey Non-additivity Type Test for Time Series Nonlinearity", %\emph{Biometrika}, \emph{72(1)}, 39-44.

\bibitem{}  Krager H., Kluger P. (1993). Nonlinearities in foreign exchange markets: a different perspective", \emph{Journal of International Money and Finance}, \emph{12}, 195-208.

\bibitem{} Lee T.H., White H., Granger C.W.J. (1993). Testing for neglected nonlinearity in time series models", \emph{Journal of Econometrics}, \emph{56}, 269–290.

\bibitem{} Lo\`{e}ve M. (1977). \emph{Probability theory 1}, 4th ed., Springer Verlag.

%\bibitem{} Lundbergh S., Terasvirta T. (2002): "Forecasting with Smooth Transition Autoregressive Models", in M.P. %Clements and D.F. Henrdry (eds.), \emph{A Comparison to Economic Forecasting}, Blackwell Publishers.

%\bibitem{} Mincer J., Zarnowitz V. (1969): "The Evaluation of Econmic Forecasts", in \emph{Economic Forecasts and Expectations}, ed. by J. Mincer, National Bureau in Economic Research, New York.

\bibitem{}  Peel D.A., Speight A.E.H. (1994). Testing for nonlinear dependence in inter-war exchange rates, \emph{Weltwiirtschaftliches Archiv}, \emph{130}, 391-417

\bibitem{} Potscher B.M., Prucha I.R. (1989). A Uniform Law of large Numbers for Dependent and Heterogeneous Data Processes,
\emph{Econometrica}, \emph{57(3)}, 675-683.

\bibitem{} Potter S. (1995). A nonlinear approach to U.S. GNP, \emph{Journal of Applied Econometrics}, \emph{10}, 109–125.

\bibitem{} Priestley, M. B. (1980). State-Dependent Models: A General Approach to Nonlinear Time Series Analysis, \emph{Journal of Time Series Analysis}, \emph{1} ,47- 71.

%\bibitem{} Skalin J., Teraesvirta T. (2000): "Modelling Asymmetries and Moving Equilibria in Unemployed Rates", \emph{Macroeconomics Dynamics}, \emph{6(2)}, 202–241.

%\bibitem{} Teraesvirta T. (1994): "Specification, Estimation and Evaluation of Smooth Transition Autoregressive Models", \emph{Journal of the American Statistical Association}, \emph{89}, 208–218.

%\bibitem{} Terasvirta T., Anderson H.M. (1992): "Characterizing Nonlinearities in Business Cycles with Smooth Transition Autoregressive Models", \emph{Journal of Applied Econometrics}, \emph{7}, 119-136.

\bibitem{} Terasvirta T., Lin C.F., Granger C.W.J. (1993). Power of the Neural Network Linearity Test, \emph{Journal of Time Series Analysis}, \emph{14}, 209–220.

\bibitem{} Tiao G.C., Tsay R.S. (1994). Some advances in nonlinear and adaptive modelling time series, \emph{Journal of Forecasting}, \emph{13}, 109-131.

\bibitem{} Tong H. (1978). On a threshold model, in Chen C.H. (ed.), \emph{Pattern Recognition and Signal Processing}, 101-141. Amsterdam: Sijhoff and Noordoff.

    \bibitem{}  Tong H., Lim K.S. (1980). Threshold autoregression, limit cycles and cyclical data, \emph{Journal of the Royal Statistical Society}, \emph{B 42}, 245-292.

\bibitem{} Tong H. (1983).Threshold Models in Nonlinear Time Series Analysis, New York, Springer – Verlag.

\bibitem{}  Tong H. (1986). On estimating thresholds in autoregressive models, \emph{Journal of Time Series Analysis}, \emph{7}, 178-190.

    \bibitem{}  Tong H. (1995). \emph{Nonlinear Time Series: A Dynamical System Approach}, Oxford University Press.

\bibitem{}  Tong H., Yeung I. (1991). On tests for Self-Exciting Threshold Autoregressive-Type nolinearity in partially observed time series, \emph{Applied Statistics}, \emph{40}, 43-62.

\bibitem{}  Tsay R.S. (1986). Nonlinearity Tests For Time Series, \emph{Biometrika}, \emph{73}, 461-466.

%\bibitem{}  Tsay R.S. (1989): "Testing and Modelling Threshold Autoregressive Process", \emph{Journal of American Statistical Association}, \emph{c. 84}, \emph{s. 405}, 231-240.

%\bibitem{}  van Dijk D., Terasvirta T., Franses P.H. (2002): "Smooth Transition Autoregressive Models - A Survey and Recent Development", \emph{Econometric Reviews}, \emph{21(1)}, 1-47.


\bibitem{} White H. (1984). \emph{Asymptotic theory for econometricians}, Academic Press, New York.

\bibitem{} White H. (1994). \emph{Estimation, Inference and Specification Analysis},
Econometric Society Monographs n. 22, Cambridge University Press:
Cambridge.

\end{thebibliography}
\end{document}